\input amstex
\documentstyle{amsppt}
\loadmsbm

\def\card#1{\vert#1\,\vert}

\def\R{\Bbb R}

\def\frac#1#2{{{#1} \over {#2}}}

\def\S{{\frak S}}

\loadbold

\def\R{{\Bbb R}}
\def\Z{{\Bbb Z}}

\def\vi{{\tilde \bold v}^{(i)}}
\def\gammai{{\gamma^{(i)}}}

\def\bfw{{\bold w}}
\def\bfv{{\bold v}}
\def\bfe{{\bold e}}

\def\frac#1#2{{{#1} \over {#2}}}

\def\card#1{{\left\vert{#1}\right\vert}}

\def\wik{{\tilde \bold w}^{(i)}(E_k)}
\def\wzj{{\tilde \bold w}^{(0)}(E_j)}
\def\wzk{{\tilde \bold w}^{(0)}(E_k)}

\def\tij{t^{(i)}(E_j)}
\def\tik{t^{(i)}(E_k)}
\def\tzj{t^{(0)}(E_j)}
\def\tzk{t^{(0)}(E_k)}

\topmatter

\title
Sequentiality Restrictions in Special Relativity
\endtitle

\rightheadtext{Sequentiality}

\author
Mark I. Heiligman
\endauthor

\address
Intelligence Advanced Research Activity (IARPA),
Office of the Director of National Intelligence,
Washington, D.C.
\endaddress

\email
mark.i.heiligman\@ugov.gov 
\endemail

\date
February 4, 2010
\enddate


\abstract 
Observers in different inertial frames can see a set of spacelike
separated events as occurring in different orders.  Various
restrictions are studied on the possible orderings of events that can
be observed.  In $1+1$-dimensional spacetime $\{(1 2 3), (2 3 1), (3 1
2)\}$ is a disallowed set of permutations. In $3+1$-dimensional
spacetime, any four different permutations on the ordering of $n$
events can be seen by four different observers, and there is a
set of four events such that any of the $4!=24$ possible orderings can
be observed in some inertial reference frame.  A more complicated
problem is that of five observers and five events, where of the
7,940,751 choices of five distinct elements from $\S_5$ (containing
the identity), all but at most one set of permutations can be
realized, and it is shown that this remaining case is impossible.  For
six events and five observers, it is shown that there are at least
7970 cases that are unrealizable, of which at least 294 do not come
from the forbidden configuration of five events.
\endabstract

\thanks
{\it Disclaimer}. 
This paper was produced in the author's personal capacity and 
all statements of fact, opinion, or analysis expressed in this paper
are solely those of the author and do not necessarily reflect the
official positions or views of the Office of the Director of National
Intelligence (ODNI) or any other government agency.  Nothing in the
content should be construed as asserting or implying U.S. Government
authentication of information or ODNI endorsement of the author's
views.
\endthanks

\endtopmatter

\document

\head
Introduction
\endhead

One of the well known ``paradoxes'' of special relativity is that two
observers who are moving in different inertial reference frames will
disagree about whether two events in space-time are simultaneous.
Consequently, they can also disagree about the sequential order in
which events in space-time occur.  The purpose of this note is to
begin investigating whether there are restrictions on the relative
orderings of these events.

Generally for $n$ events in spacetime, an observer will see them
ocurring in some order, which can be viewed as an element of $\S_n$,
the symmetric group of permutations of the set $\{1,2,\ldots,n\}$.
Each intertial frame of reference will assign an element $\pi\in
\S_n$ to the observed sequence of events.  The basic question is
whether there are restrictions on which subsets of $\S_n$ can occur for
multiple observers of the same $n$ events.

In all that follows, the speed of light $c$ will be taken to be 1, and all
velocities will be taken as the fraction of the speed of light.

\subhead
Historical Note
\endsubhead

A number of these ideas were discussed in Professor Richard Stanley's
paper \cite{S}.  In that paper, he mentions a mathematician from NSA
whose name he can't recall who made the connection between event
orderings in Minkowski space with hyperplane arrangements.  The author
of the present paper is that mathematician, who at the time of the
discussion was working at NSA. The current author truly appreciates
the reference from \cite{S}, and was deeply gratified in the course of
checking out the background references for this paper to find that
Professor Stanley found the question to be of sufficient interest
to write on the topic.

\head
The 1-Dimensional Case
\endhead

Consider the situation of events on a line, so that their space-time
coordinates consist only of an $x$-coordinate and a time. For an observer
in an inertial frame, the $n$ events have coordinates $(x_i, t_i)$ for 
$i=1,2,\ldots,n$. Any other inertial frame is related to this rest frame
by a characteristic velocity $v$, in the range $-1 < v < 1$.

For another observer in an inertial frame moving
at velocity $v$ with respect to the initial frame, the $i$-th event will
have coordinates $(x_i^\prime, t_i^\prime)$ where
$$\eqalign{
x_i^\prime&=\gamma\,(x_i-v\,t_i) \cr
t_i^\prime&=\gamma\,(t_i-v\,x_i) \cr
}$$
where
$$\gamma = \frac{1}{\sqrt{1-v^2}} $$
is the Lorentz factor.

Now suppose that $\sigma\in \S_n$ is a permutation, and set
$t_i=i$ and $x_i=2\,\bigl(i-\sigma(i)\bigr)$ and let $v=\frac{1}{2}$.
Then $t_i^\prime = \gamma\,\sigma(i)$, so while an observer in the
rest frame sees the events in order $1,2,\ldots,n$, an observer in an
inertial frame moving with velocity $v=\frac{1}{2}$ with respect to
the rest frame sees the events in order
$\sigma(1),\sigma(2),\ldots,\sigma(n)$. This shows that for any pair
of permutations in $\S_n$ there are $n$ events and a pair of inertial
frames such that each observer sees the events with their respective
permutations on the order of occurrence.  

On the other hand, consider the case of three events with three
observers and whether it is possible for one observer to see the
events in order $(1,2,3)$, the second observer to see the events in
order $(2,3,1)$ and the third observer to see the events in order
$(3,1,2)$. From the perspective of at any observer who numbers the
events as $(1,2,3)$, the other two observers number the events as
$(2,3,1)$ and $(3,1,2)$.  Also, for one of the observers, which will
be the rest frame for the purposes of this discussion, the other two
reference frames are moving with velocities $v^{(1)}$ and $v^{(2)}$
that are both positive. If the events have coordinates $(x_i, t_i)$ for 
$i=1,2,3$ in the rest frame, then $t_1<t_2<t_3$, while in the other
two frames
$$\eqalign{
t_2^{(1)}=\gamma^{(1)}\bigl(t_2-v^{(1)}\,x_2\bigr)
&<t_3^{(1)}=\gamma^{(1)}\bigl(t_3-v^{(1)}\,x_3\bigr)
<t_1^{(1)}=\gamma^{(1)}\bigl(t_1-v^{(1)}\,x_1\bigr) \cr
t_3^{(2)}=\gamma^{(2)}\bigl(t_3-v^{(2)}\,x_3\bigr)
&<t_1^{(2)}=\gamma^{(2)}\bigl(t_1-v^{(2)}\,x_1\bigr)
<t_2^{(2)}=\gamma^{(2)}\bigl(t_2-v^{(2)}\,x_2\bigr) \cr
}$$
or simply extracting out some relevant inequalities
$$\eqalign{
t_2-v^{(1)}\,x_2 &< t_1-v^{(1)}\,x_1 \cr
t_3-v^{(1)}\,x_3 &< t_1-v^{(1)}\,x_1 \cr
t_3-v^{(2)}\,x_3 &< t_1-v^{(2)}\,x_1 \cr
t_3-v^{(2)}\,x_3 &< t_2-v^{(2)}\,x_2. \cr
}$$
Now adding the inequlaties $t_1<t_2<t_3$ appropriately gives
$$\eqalign{
t_1+t_2-v^{(1)}\,x_2 &< t_2+t_1-v^{(1)}\,x_1 \cr
t_1+t_3-v^{(1)}\,x_3 &< t_3+t_1-v^{(1)}\,x_1 \cr
t_1+t_3-v^{(2)}\,x_3 &< t_3+t_1-v^{(2)}\,x_1 \cr
t_2+t_3-v^{(2)}\,x_3 &< t_3+t_2-v^{(2)}\,x_2. \cr
}$$
so that
$$\eqalign{
v^{(1)}\,x_1 &< v^{(1)}\,x_2 \cr
v^{(1)}\,x_1 &< v^{(1)}\,x_3 \cr
v^{(2)}\,x_1 &< v^{(2)}\,x_3 \cr
v^{(2)}\,x_2 &< v^{(2)}\,x_3 \cr
}$$
and therefore $x_1 < x_2 < x_3$ since $v^{(1)}>0$ and $v^{(2)}>0$.

Now notice that 
$\gamma^{(1)}\,\bigl(t_2-v^{(1)}\,x_2\bigr) 
< \gamma^{(1)}\,\bigl(t_3-v^{(1)}\,x_3\bigr)$
and
$\gamma^{(2)}\,\bigl(t_3-v^{(2)}\,x_3\bigr) 
< \gamma^{(2)}\,\bigl(t_2-v^{(2)}\,x_2\bigr)$
implies
$$
\bigl(t_2-v^{(1)}\,x_2\bigr) + \bigl(t_3-v^{(2)}\,x_3\bigr)
<\bigl(t_3-v^{(1)}\,x_3\bigr) + \bigl(t_2-v^{(2)}\,x_2\bigr)
$$
and therefore
$$\bigl(v^{(1)} - v^{(2)}\bigr)\,\bigl(x_3 - x_2\bigr)
=v^{(1)}\,x_3+v^{(2)}\,x_2
-v^{(1)}\,x_2 - v^{(2)}\,x_3
<0$$
which implies $v^{(1)} - v^{(2)} < 0$ since $x_3 - x_2 > 0$.

On the other hand, 
$\gamma^{(1)}\bigl(t_2-v^{(1)}\,x_2\bigr) < \gamma^{(1)}\bigl(t_1-v^{(1)}\,x_1\bigr)$
and
$\gamma^{(2)}\bigl(t_1-v^{(2)}\,x_1\bigr) < \gamma^{(2)}\bigl(t_2-v^{(2)}\,x_2\bigr)$
implies
$$
\bigl(t_2-v^{(1)}\,x_2\bigr) + \bigl(t_1-v^{(2)}\,x_1\bigr)
<\bigl(t_1-v^{(1)}\,x_1\bigr) + \bigl(t_2-v^{(2)}\,x_2\bigr)
$$
and therefore
$$\bigl(v^{(2)} - v^{(1)}\bigr)\,\bigl(x_2 - x_1\bigr)
=v^{(2)}\,x_2 + v^{(1)}\,x_1 - v^{(1)}\,x_2 - v^{(2)}\,x_1
< 0
$$
which implies $v^{(2)} - v^{(1)} < 0$ since since $x_2 - x_1 > 0$.

This gives a contradiction, thereby showing that the triple of permutations
$\{(1,2,3),(2,3,1),(3,1,2)\}$ is not realizable by three observers in $1+1$
relativistic spacetime.

It might be guessed that this example provides an essential obstruction,
and indeed any set of three permutations that contains a subset of 
events that follows the pattern $\{(1,2,3),(2,3,1),(3,1,2)\}$ is
unrealizable. However, it is not the case that all such unrealizable sets
of three permutations contain a subset of this type.

Notice that the essential issue is that as the velocity increases, if
the order of a pair of events flips, then it can't flip back.  For
example, the three permutations on six events
$\{(1,2,3,4,5,6),(2,1,3,4,6,5),(1,2,4,3,6,5)\}$ is a set of
permutations that cannot occur in 1+1 Minkowski spacetime because
regardles of how they are ordered, one of the pairs $\{1,2\}$ or
$\{3,4\}$ or $\{5,6\}$ flips its order and then flips back.  This 
example is interesting because there is no subset of the events
that follows the pattern $\{(1,2,3),(2,3,1),(3,1,2)\}$.

If $\pi\in\S_n$, a reversal within $\pi$ is a pair $(i,j)$ with $i<j$ such that
$\pi(j)<\pi(i)$.  Let $\rho(\pi)$ denote the set of reversals, i.e.
$$\rho(\pi) = \{(i,j)\,\vert\, 1\le i < j \le n \,\,\&\,\,\pi(j) <
\pi(i)\}.$$ 
If $\pi_1,\pi_2\in\S_n$, write $\pi_1<\pi_2$ if
$\rho(\pi_1)\subset\rho(\pi_2)$.  This gives a partial order on
$\S_n$. With respect to this ordering, the
minimal element is the identity permutation $\pi^{(0)}=(1,2,\ldots,n-1,n)$ and
the maximal element is the reverse of the identity $(n,n-1,\ldots,2,1)$.

Without any loss of generality, in considering whether a set
of permutations is realizable, it may be assumed that one of the
permutations is the identity.
It is easily seen that a necessary condition for a set
$\{\pi_1,\pi_2,\ldots,\pi_k\}$ of permutations to be realizable
is that they can be put in order with $\pi^{(0)}=\pi_1<\pi_2,\ldots<\pi_k$.
This was pointed out in \cite{S}.

It might be hoped that this criterion, namely that there is an
ordering of permutations, would also be a sufficient condition for a
set of permutations to be realizable.  Alas, as pointed out in
\cite{S}, the paper \cite{GP} gives a counterexample of a set of
five permutations $\{\pi_0,\pi_1,\pi_2,\pi_3,\pi_4\}$ with $\pi_0$ the
identity and with
$\rho(\pi_1)\subset\rho(\pi_2)\subset\rho(\pi_3)\subset\rho(\pi_3)$
that is not realizable.

Nevertheless, an interesting question and perhaps one that is easier,
would be to consider a set of three observers in $1+1$-dimensional
Minkowski space, and assume that one of the observers in the rest
frame sees the events in order $\pi^{(0)}$ and that the other two
observers see the events in order $\pi_1$ and $\pi_2$, and that
$\pi_1<\pi_2$. Then is it possible to find an arrangement of $n$
events in spacetime and a pair of velocities for observers 1 and 2
that realize these two permutations? The example from \cite{GP}
is not a counterexample to this question, since it involves more than
three observers.

\head
The General Case
\endhead

Consider the situation of $m+1$ observers of $n$ events. The $i$-th
observer sees the $k$-th event $E_k$ as having space-time coordinates
$\bigl(\tik,\wik\bigr)$.  Assume that the observers are numbered
$i=0,1,\ldots,m$ and that the events are numbered $j=1,\ldots,n$.
If $\vi$ is the velocity of observer $i$ as seen by observer 0, then
the Lorentz transformations give
$$\eqalign{
\tik &= \gammai\,\bigl(\tzk - \vi\cdot\wzk\bigr) \cr
}$$
where
$$\gammai = \gamma(\vi) = {1 \over \sqrt{1-\vi\cdot\vi}}$$
is the Lorentz contraction factor.  In all of this,  
$\vi\cdot\vi < 1$, with the
speed of the $i$-th observer being $\sqrt{\vi\cdot\vi}$ as
a fraction of the speed of light.  Note also that $\gammai > 0$ for
all $i$.

If $\tzj < \tzk$, but $\tzk - \vi\cdot\wzk < \tzj - \vi\cdot\wzj$,
then observer $i$ sees event $E_k$ as preceeding event $E_j$ while the
observer in the rest frame sees event $E_j$ as preceeding event $E_k$.
The equation $\tzk - \tzj = \vi\cdot\bigl(\wzk - \wzj\bigr)$
determines a hyperplane in the space of velocities, which separate
observers into two sets depending on whether $E_j$ is before or after
$E_k$. If there are $n$ points in spacetime, then there are ${n
\choose 2}$ such separating hyperplanes. Each region of such a hyperplane 
arrangment corresponds to a different set of time orderings of the
events to the different observers, i.e. to different subsets of $\S_n$.

To establish some basic notation, for an observer in reference frame
$j$, the time of an event $E$ is $t^{(j)}(E)$ while the spatial
coordinates of event $E$ are $\tilde \bfw^{(j)}(E)$.  From the point of view
of reference frame 0, the velocity of reference frame $j$ is $\tilde \bfv^{(j)}$
and the corresponding Lorentz factor is
$$\gamma^{(j)} = {1 \over \sqrt{1-\tilde \bfv^{(j)}\cdot\tilde \bfv^{(j)}}}$$
which is always positive.  The Lorentz transformation then gives
$$t^{(j)}(E) = \gamma^{(j)}\left(t^{(0)}(E)-\tilde \bfv^{(j)}\cdot\tilde \bfw^{(0)}(E)\right).$$

\definition{Definition}
If $t^{(j)}(E_{i_1}) < t^{(j)}(E_{i_2}) < \cdots < t^{(j)}(E_{i_n})$,
then the reference frame $j$, or the observer $j$, sees the events in
order $(i_1, i_2,\ldots, i_n)\in\S_n$
\enddefinition

Each observer sees the $n$ events under
consideration in some order that can be described by an element
$\sigma\in \S_n$.  For the $i$-th observer, suppose that the order in
which the events are observered is given by the permutation
$\sigma_i\in \S_n$, so that
$$t^{(i)}(E_{\sigma_i(1)}) < t^{(i)}(E_{\sigma_i(2)}) < \ldots 
< t^{(i)}(E_{\sigma_i(n)})$$
which can be expressed as saying that if $j<k$ then
$t^{(i)}(E_{\sigma_i(j)}) < t^{(i)}(E_{\sigma_i(k)})$.
Without loss of generality, it may be assumed that the events are
labeled so that observer 0 sees them in order $1,2,\ldots,n$ so that
that $t^{(0)}(E_{1}) < t^{(0)}(E_{2}) < \ldots < t^{(0)}(E_{n})$
which means that $\sigma_0\in \S_n$ is the identity permutation. 

\definition{Definition}
A set $Q=\{\pi_1,\ldots,\pi_k\}\subset\S_n$ of cardinality $k$ is said
to be {\it realizable} in $d+1$-dimensional Minkowski space if there
is a set of $n$ events in $d+1$-dimensional Minkowski space and a set
of $k$ inertial reference frames in $d+1$-dimensional Minkowski space
such that observer $i$ sees the $n$ events in order $\pi_i$ for
$i=1,\ldots,k$.
\enddefinition

It is important to realize that there is no preferred labeling of the events,
so if the events are relabeled by $\sigma^{-1}\in\S_n$, the set $Q$ becomes
$Q\,\sigma = \{\pi_1\,\sigma,\ldots,\pi_k\,\sigma\}$.

\definition{Definition}
If $Q\subset\S_n$ and if $\sigma\in\S_n$ then the set $Q\,\sigma$ is
said to be equivalent to $Q$.  The equivalent sets to $Q$ that contain
the identity permutation are of the form $Q\,\pi^{-1}$ as $\pi$ ranges
over the elements of $Q$. If $Q$ is actually a subgroup of $\S_n$ then
the only subset of $\S_n$ that contains the identity and is equivalent
to $Q$ is $Q$ itself.
\enddefinition

Since special relativity is invariant under time reversal, it is
useful to consider the time reversal of a set of events. Let
$\pi_r=(n, n-1,\ldots,2,1)$ be the permutation that reverses the labeling
of all the elements.  Applying $\pi_r$ on the left simply reverses the order
of any permutation.

\definition{Definitions}
If $Q\subset\S_n$, its {\it time reversal} is the set $\pi_r\,Q\pi_r$. 
The set $Q$ is {\it time reversal invariant} if $\pi_r\,Q\pi_r=Q$.
\enddefinition

\subhead
Counting invariant sets
\endsubhead

It is intersting to count the number of time reversal invariant sets
of $\S_n$ of size $k$ that contain the identity permutation
$\pi_0$. 

The permutation $\pi_r$ is an involution (i.e. $\pi_r^2=1$).
Let $C$ be the set of elements of $\S_n$ that are invariant under
conjugation by $\pi_r$ and let $c=\card{C}$ be the number of elements
of $\S_n$ that are invariant under conjugation by $\pi_r$. The number
of pairs of elements of $\S_n$ that are conjugates under $\pi_r$ is
then $c^\prime = (n!-c)/2$ and $c^\prime=\card{C^\prime}$ where
$C^\prime$ is the set of pairs of elements of $\S_n$ that are conjugate
under $\pi_r$.

If $\pi_r\,Q\pi_r=Q$, then the elements of $Q$ are either invariant
under conjugation by $\pi_r$ or come in pairs that are conjugates
under $\pi_r$.  Let $i$ be the number of elements of $Q$ that are
invariant under conjugation by $\pi_r$, so $j=(k-i)/2$ is the number of
conjugate pairs in $Q$. The number of ways of choosing $i$ elements of
$C$, where one of the elements is $\pi_0$, the identity permutation is
$c-1 \choose i-1$ and the number of ways of choosing $j=(k-i)/2$
elements of $C^\prime$ is $c^\prime \choose j$. Since $i=k-2j$, the total
number of time reversal invariant sets $Q$ is 
$$\sum_{j=0}^J {c^\prime \choose j}\,{c-1 \choose k-2j-1}$$
where $J=(k-1)/2$ if $k$ is odd, and $J=(k-2)/2$ if $k$ is even, and in
either case $J=\bigl[\frac{k-1}{2}\bigr]$.
Since $c^\prime = (n!-c)/2$, all that remains is to find $c$.

From the group theory perspective, $C = C(\pi_r)$ is the centralizer
of $\pi_r$, which is the subgroup of all permutations that commute
with $\pi_r$, so that $c=\card{C(\pi_r)}$. The right cosets of
$C(\pi_r)$ are in one-to-one correspondence with conjugacy classes of
$\pi_r$ in $\S_n$.  In turn, the conjugacy class of $\pi_r$ in $\S_n$
consists of all permutations having the same cycle structure as
$\pi_r$.

If $n$ is even then $\pi_r$ is a product of $\frac{n}{2}$ 2-cycles,
i.e. $\pi_r=(1\,\, n)\,(2\,\, n-1)\cdots(\frac{n}{2}\,
\frac{n}{2}+1)$, while if $n$ is odd then $\pi_r$ is a product of
$\frac{n-1}{2}$ 2-cycles and a single 1-cycle, i.e.  $\pi_r=(1\,\,
n)\,(2\,\, n-1)\cdots(\frac{n-1}{2}\,
\frac{n+3}{2})\,(\frac{n+1}{2})$.

To determine the number of ways that $\bigl[\frac{n}{2}\bigr]$ 2-cyles
can be made, note that there are ${n \choose 2}$ ways of picking the
first 2-cycle, and then after that there are ${n-2 \choose 2}$ ways of
picking the second 2-cycle, after which there are ${n-4 \choose 2}$ ways
of picking the third 2-cyle, etc. for a total of 
$${n \choose 2}\cdot{n-2 \choose 2}\cdot{n-4 \choose 2}
  \cdots{n-2[\frac{n-2}{2}] \choose 2}
=\frac{n!}{2^{[\frac{n}{2}]}}$$
ways of picking the needed numbers of 2-cycles.  However, this needs to be
divided by the number of ways that these two cyles can be ordered, since
such a reordering gives rise to the same permutation in the end. Thus the
total number of permutations having the same cycle structure as $\pi_r$ is
$n!/\bigl([\frac{n}{2}]!\,2^{[\frac{n}{2}]}\bigr)$,
which is the size of the conjugacy class of $\pi_r$ in $\S_n$. Therefore
the size of the centralizer of $\pi_r$ is
$$\hbox{$c=\bigl[\frac{n}{2}\bigr]!\,2^{[\frac{n}{2}]}.$}$$

For $n=5$, $c=8$ and $c^\prime = 56$, while for $n=6$, $c=48$ and $c^\prime=336$.
For $k=5$ observers $J=2$ and therefore for $n=5$, the number of time reversal invariant
subsets of $\S_5$ of size 5 is 
$${48 \choose 0}{7 \choose 4}+{48 \choose 1}{7 \choose 2}+{48 \choose 2}{7 \choose 0}
=2171
$$
while for $n=6$, the number of time reversal invariant
subsets of $\S_6$ of size 5 is
$${336 \choose 0}{47 \choose 4}+{336 \choose 1}{47 \choose 2}+{336 \choose 2}{47 \choose 0}
=597861
$$

\subhead
Linearly dependent velocities
\endsubhead

Now suppose that there is a linear dependence between the velocities
of all the observers, so that
$$\sum_{i=1}^{m} \alpha_i\,\vi = 0$$
for some real numbers $\alpha_i\in\R$. (Note that ${\tilde \bold v}^{(0)} = 0$ so
that $\alpha_0$ is irrelevant in this linear relation.)  Then
$$\eqalign{
\sum_{i=1}^{m} \gammai^{-1}\,\alpha_i\,\tik
&=\sum_{i=1}^{m} \alpha_i\,\bigl(\tzk - \vi\cdot\wzk\bigr) \cr
&=\left(\sum_{i=1}^{m} \alpha_i\right)\,\tzk 
 - \left(\sum_{i=1}^{m} \alpha_i\,\vi\right)\cdot\wzk \cr
&=A\,\,\tzk \cr
}$$

It is important to assume now that $\alpha_i \ne 0$ for $i=1,\ldots,m$ 
and that $A \ne 0$, as well.  These actually are quite mild assumptions,
since the set of $m$ velocities that will satisfy any set of strict sequentiality
inequalities will be an open set, slight variations in the velocities will
not affect the ordering of events for the various observers.  These slight
variations in the $\vi$'s will be enough to insure that without any loss of
generality, it should be possible to take all the $\alpha_i$'s to be non-zero
and at the same time also take their sum $A$ to be nonzero.

With this assumption in hand, it is also clear that the $\alpha_i$'s
can all be multiplied by any nonzero real number, and the linear
dependence of the velocities will remain.  Without loss of generality,
then, it can be assumed that the $\alpha_i$'s are such that $A>0$.
Note that this means that the $\alpha_i$'s cannot all be negative.
This is the first restriction on the signs of the $\alpha_i$'s.

It follows that $\tzj<\tzk$ for $j<k$ and therefore
$$0 \le A\,\left(\tzk-\tzj\right) = 
\sum_{i=1}^{m} \gammai^{-1}\,\alpha_i\,\left(\tik - \tij\right)$$
for $j < k$.

For a given $j$ and $k$ with $1 \le j < k \le n$, let
$$\eqalign{
I(j,k) &= \{i \mid \sigma_i^{-1}(j) > \sigma_i^{-1}(k)\} \cr
\bar I(j,k) &= \{i \mid \sigma_i^{-1}(j) < \sigma_i^{-1}(k)\} \cr
}$$
and consider the following sign pattern of the $\alpha_i$'s:
$$\eqalign{
\alpha_i < 0 &\qquad\hbox{if and only if}\qquad i\in I(j,k) \cr
\alpha_i > 0 &\qquad\hbox{if and only if}\qquad i\in \bar I(j,k) \cr
}$$

Since $j<k$ if and only if $t^{(i)}(E_{\sigma_i(j)}) <
t^{(i)}(E_{\sigma_i(k)})$, it also follows that $\tzj < \tzk$ if and only
if $\sigma_i^{-1}(j) < \sigma_i^{-1}(k)$, and therefore $\tik - \tij >
0$ if and only if $i\in I(j,k)$.  Consequently
$\alpha_i\bigl(\tik-\tij\bigr)<0$ for all $i\in I(j,k)$.  Similarly,
$\tzj > \tzk$ if and only if $\sigma_i^{-1}(j) > \sigma_i^{-1}(k)$,
and therefore $\tik - \tij < 0$ if and only if $i\in \bar I(j,k)$.
Consequently $\alpha_i\bigl(\tik-\tij\bigr)<0$ for all $i\in \bar
I(j,k)$.  Therefore $\alpha_i\bigl(\tik-\tij\bigr)<0$ for all $i$
pairs, unconditionally.

As a result of this, the above sign pattern leads to the contradiction
$$0 \le \sum_{i=1}^{m} \gammai^{-1}\,\alpha_i\,\left(\tik - \tij\right) 
< 0$$
since the $\gammai$'s are all positive.  The consequence is that each
$(j,k)$ pair with $j<k$ invalidates one possible sign pattern for the
$\alpha_i$'s.

The total number of sign patterns up for consideration is $2^{m-1}-1$
since the all negative sign pattern has already been eliminated by the
requirement $A>0$.  The number of $(j,k)$ pairs is $n \choose 2$, so
if ${n \choose 2} \ge 2^{m-1}-1$, it is possible that all sign patterns
are eliminated.  The conclusion is that in such a case, the suggested
sequence of observations of the events by the different observers is not
possible.

It is interesting to consider the solutions of the Diophantine equation
${n \choose 2} = 2^{m-1}-1$.  The known solutions are $\{(n=3,m=3),(n=6,m=5),
(n=91,m=13)\}$.  The first corresponds to restrictions on one space dimension,
which has already been discussed.  The second solution corresponds to 
restrictions on $3+1$-dimensional space-time, while the third solution would
correspond to restrictions on $11+1$-dimensional space-time.  Of course there
could be other restictions, since all that is necessary is that the inequality
be satisfied.

\subhead
Four Events in 3-Dimensional Space
\endsubhead

If there are four events in spacetime, labeled $E_1,E_2,E_3,E_4$ an
observer will see the sequencing of these events in one of 24 possible
ways.  The goal of this section is to show that in 3-space, it is possible
to choose these four events once and for all so they have a 
fixed sequence $t^{(0)}(E_1) < t^{(0)}(E_2) < t^{(0)}(E_3) < t^{(0)}(E_4)$ in reference
frame 0, and that for any $\sigma\in \S_4$, there is an inertial frame that 
sees these events in order $\sigma$, i.e. that 
$$t^{(\sigma)}(E_{\sigma(1)}) < t^{(\sigma)}(E_{\sigma(2)}) 
< t^{(\sigma)}(E_{\sigma(3)}) < t^{(\sigma)}(E_{\sigma(4)}).$$ 
The velocity of this reference frame with respect
to reference frame 1 is $\tilde \bfv^{(\sigma)}$ and
$$t^{(\sigma)}(E_i) 
= \gamma^{(\sigma)}\,\left(t^{(0)}(E_i)-\tilde \bfv^{(\sigma)}\cdot\tilde \bfw^{(0)}(E_i)\right)$$
for $i=1,2,3,4$.  For this purpose, it is possible to take
$$\eqalign{
\tilde \bfw^{(0)}(E_1)=(0,0,0) &\qquad t^{(0)}(E_1)=1 \cr
\tilde \bfw^{(0)}(E_2)=(1,0,0) &\qquad t^{(0)}(E_2)=2 \cr
\tilde \bfw^{(0)}(E_3)=(0,1,0) &\qquad t^{(0)}(E_3)=3 \cr
\tilde \bfw^{(0)}(E_4)=(0,0,1) &\qquad t^{(0)}(E_4)=4 \cr
}$$
i.e. $t^{(0)}(E_i)=i$ for $i=1,2,3,4$.  The Lorentz transformations then give
$$t^{(\sigma)}(E_i) 
= \gamma^{(\sigma)}\left(i-\tilde \bfv^{(\sigma)}\cdot\tilde \bfw^{(0)}(E_i)\right)$$
so
$$t^{(\sigma)}(E_{\sigma(i)}) 
= \gamma^{(\sigma)}\left(\sigma(i)
-\tilde \bfv^{(\sigma)}\cdot\tilde \bfw^{(0)}(E_{\sigma(i)})\right)$$
for $i=1,2,3,4$.  
Since $\gamma^{(\sigma)}>0$, the conditions that
$t^{(\sigma)}(E_{\sigma(i)}) < t^{(\sigma)}(E_{\sigma(i+1)})$  be met
for $i=1,2,3$ can be achieved by setting 
$$t^{(\sigma)}(E_{i}) = \gamma^{(\sigma)}\,(s_\sigma+\sigma^{-1}(i))$$
for some constant $s_\sigma$ to be determined that depends only on
$\sigma\in \S_4$, but not on $i$.  Then
$$t^{(\sigma)}(E_{\sigma(i)}) = \gamma^{(\sigma)}(s_\sigma+i)$$
and the desired ordering follows since
$$\gamma^{(\sigma)}\,(s_\sigma+1) 
<\gamma^{(\sigma)}\,(s_\sigma+2)
<\gamma^{(\sigma)}\,(s_\sigma+3)
<\gamma^{(\sigma)}\,(s_\sigma+4).$$
It follows that the equations to be satisfied are
$$i-\tilde \bfv^{(\sigma)}\cdot\tilde \bfw^{(0)}(E_i)
=s_\sigma+\sigma^{-1}(i)$$
for $i=1,2,3,4$.  Since $\tilde \bfw^{(0)}(E_1) = 0$ it follows that
$$s_\sigma = 1-\sigma^{-1}(1)$$
so that the remaining equations to be solved are
$$\tilde \bfv^{(\sigma)}\cdot\tilde \bfw^{(0)}(E_i) = i-1+\sigma^{-1}(1)-\sigma^{-1}(i)$$
for $i=2,3,4$, where the unknowns are the three components of
$\tilde \bfv^{(\sigma)}$.  However, this follows immediately from
noting that $\tilde \bfw^{(0)}(E_2),\tilde \bfw^{(0)}(E_3),\tilde \bfw^{(0)}(E_4)$ are just the
basis vectors $(1,0,0),(0,1,0),(0,0,1)$, respectively.  This gives
that the desired velocity for an observer to see the permutation
$\sigma$ on the temporal ordering of events $E_1,E_2,E_3,E_4$ is
$$\tilde \bfv^{(\sigma)} 
= \bigl(1+\sigma^{-1}(1)-\sigma^{-1}(2),
2+\sigma^{-1}(1)-\sigma^{-1}(3),
3+\sigma^{-1}(1)-\sigma^{-1}(4)\bigr).$$

A slight problem here is that $\card{\tilde \bfv^{(\sigma)}}$ can be bigger
than 1.  In fact its maximum value is achieved when $\sigma^{-1}(1)=4$,
$\sigma^{-1}(2)=3$, $\sigma^{-1}(3)=2$, $\sigma^{-1}(4)=1$, in which
case $\tilde \bfv^{(\sigma)}\cdot\tilde \bfv^{(\sigma)} = 56$.
Replacing $\tilde \bfv^{(\sigma)}$ by $\tilde \bfv^{(\sigma)}/8$ ensures that
$\tilde \bfv^{(\sigma)}\cdot\tilde \bfv^{(\sigma)} < 1$, as desired.  This is
compensated for by replacing $\tilde \bfw^{(0)}(E_i)$ for $i=1,2,3,4$ by
$8\,\tilde \bfw^{(0)}(E_i)$.

The general case here is that in $n$-dimensional space there can be
$n+1$ events in space-time such that for every $\sigma\in \S_{n+1}$
there is a reference frame such that the temporal ordering of the
events is $\sigma$.  Let $\bfe_i$ denote the vector in $\R^n$ with a
$1$ in coordinate $i$ and $0$'s everywhere else.  Now consider the
events $E_i$ where in reference frame 0, $\tilde \bfw^{(0)}(E_1)={\bold 0}$ and
$\tilde \bfw^{(0)}(E_i)=\bfe_{i-1}$ for $i=2,\ldots,n+1$ and $t_1(E_i)=i$ for
$i=1,\ldots,n+1$.  Then choosing $\tilde \bfv^{(\sigma)}$ so that its $i$-th
component is given by
$\bfv^{(\sigma)}_i \cdot \bfe_i= i+\sigma^{-1}(1)-\sigma^{-1}(i+1)$
means that an observer moving with velocity $\tilde \bfv^{(\sigma)}$ see the
events $E_i$ for $i=1,\ldots,n+1$ occuring in time order $\sigma$.
Here the maximum value of $\tilde \bfv^{(\sigma)}\cdot\tilde \bfv^{(\sigma)}$ is
achieved for $\sigma^{-1}(i)=n+2-i$ for $i=1,\ldots,n+1$ in which case
$\tilde \bfv^{(\sigma)}\cdot\tilde \bfv^{(\sigma)}=i+(n+1)-(n+1-i)=2\,i$
so that 
$$\hbox{$\tilde \bfv^{(\sigma)}\cdot\tilde \bfv^{(\sigma)} 
= \sum_{i=1}^n(2\,i)^2 
= \frac{4\,n\,(n-1)\,(2\,n-1)}{6}
< 2\,n^3.$}$$
This allows replacing $\tilde \bfw^{(0)}(E_i)$ by 
$\tilde \bfw^{(0)}(E_i)\cdot 2\,n^3$ and
$\tilde \bfv^{(\sigma)}$ by $\tilde \bfv^{(\sigma)}/(2\,n^3)$ for $i=1,\ldots,n+1$ 
in order to assure that $\tilde \bfv^{(\sigma)}\cdot\tilde \bfv^{(\sigma)} < 1$.

\proclaim{Proposition}
In $n+1$-dimensional Minkowski space there exists an ordered set of
$n+1$ events such that for every $\sigma\in\S_{n+1}$ there exists an
inertial observation frame that sees these $n+1$ events in the
sequence $\sigma$.
\endproclaim

\subhead
Four Observers in 3-Dimensional Space
\endsubhead

Consider now the question of whether an arbitrary number of events can
be seen by a limited number of observers in any order. So suppose that
there are $n$ events in 3-dimensional space-time, $E_i$ for
$i=1,\ldots,n$.  The goal of this section is to show that for any
$\sigma_0,\sigma_1,\sigma_2,\sigma_3 \in \S_n$ there is some way to assign
space and time coordinates to the $E_i$ and to choose a set of reference
frames $R_0,R_1,R_2,R_3$ so that reference frame $R_i$ sees the events
in time order $\sigma_i$ for $i=1,2,3,4$.  

Without loss of generality, it may be assumed that $\sigma_0$ is the
identity permutation, $\sigma_0(i)=i$ for $i=1,\ldots,n$.  In
reference frame 0, assume that the time coordinate of the $E_i$ is
given by $t^{(0)}(E_i)=i$.  Also assume that from the perspective of
reference frame $i$, the velocity of $R_i$ is given by $\tilde \bfv_i$, where
$$\eqalign{
\tilde \bfv^{(0)}&=(0,0,0) \cr
\tilde \bfv^{(1)}&=(v_1,0,0) \cr
\tilde \bfv^{(2)}&=(0,v_2,0) \cr
\tilde \bfv^{(3)}&=(0,0,v_3) \cr
}$$
for some $v\in\R$ with $0<v_1,v_2,v_3<1$, i.e. 
$\tilde \bfv^{(0)}={\bold 0}$ and
$\tilde \bfv^{(j)}=v_{j}\,\bfe_{j}$ for $j=1,2,3$

In order that $t^{(j)}(E_i)= \gamma^{(j)}\,\sigma_j(i)$, which would
mean that observer $j$ sees the events in time order $\sigma_j$,
it suffices that
$$\eqalign{\gamma^{(j)}\,\sigma_j(i) = t_j(E_i) 
&= \gamma^{(j)}\left(t_1(E_i)-\tilde \bfv^{(j)}\cdot\tilde \bfw^{(0)}(E_i)\right) \cr
&= \gamma^{(j)}\left(i-v_{j-1}\,\bfe^{(j-1)}\cdot\tilde \bfw^{(0)}(E_i)\right) \cr
}$$
so that
$$\bfe^{(j-1)}\cdot\tilde \bfw^{(0)}(E_i) = v_{j-1}^{-1}\,\bigl(i-\sigma_j(i)\bigr)$$
i.e.
$$\tilde \bfw^{(0)}(E_i) = 
\bigl( (i-\sigma_2(i))/v_1, (i-\sigma_3(i))/v_2, (i-\sigma_4(i))/v_3 \bigr).$$
This then gives the space coordinates of the event $E_i$ in the system of 
reference frame 1 that yield the desires space coordinates in the other
reference frames.

Clearly, this generalizes to $n$-dimensional space, where any $n+1$ observers
can see any number of events in arbitrary orders.

\proclaim{Proposition}
For any $Q\subset\S_k$ with $\card{Q}=n+1$, there exists a set of $k$ events
in $n+1$-dimensional Minkowski space and a set of $n+1$ inertial observation
frames that realize the set $Q$, i.e. for each $\sigma\in Q$ there is a frame
that sees the $k$ events in order $\sigma$.
\endproclaim

\subhead
Five Observers of Six Events in 3-Dimensional Space: A Computer Search
\endsubhead

In this case, it is at least possible to eliminate some selections
of five permutations from $\S_6$ simply on the basis of sign patterns.

The following tables contain the results of a computer search for
disallowed sets of oberservation sequences for five inertial
observers, each observing the same set of six events, in potentially
different orders. Without any loss of generality, it may be assumed
that one of the observers sees the events in serial order
$(1,2,3,4,5,6)$ so that rather than ${720 \choose 5} =
1,590,145,128,144$ sets of 5 permutations to search through, it is
only necessary to search through a set of ${719 \choose 4} =
11,042,674,501$, which is a much more tractable number. As in the
case of 5 events with 5 observers, any one of the observers can
be designated as having the identity permutation, so many of these
sets are equivalent.  The number of distinct equivalence classes of
sets is 2,208,534,929 by direct count.
Note that $5 \cdot 2208534929 - 11042674501 = 144$, which is total
number of 5-cycles in $\S_6$, i.e. the number of distinct subgroups of
$\S_6$ of order 5. (To see this note that there are 6 ways of picking
out 5 elements of $\{1,2,3,4,5,6\}$ and that each such set of 5 
elements gives rise to $4!=24$ distinct 5-cycles.)

The computer search for disallowed sets was based on eliminating all
possible sign patterns for the relative velocities of the different
observers. The total number of disallowed sign patterns found by the
computer search was 294. A closer examination of these 294 answers
shows that the actual number can be substantially reduced. If $S =
\{\pi_0,\pi_1,\pi_2,\pi_3,\pi_4\}$ is any disallowed pattern, then so
is $S\sigma =
\{\pi_0\sigma,\pi_1\sigma,\pi_2\sigma,\pi_3\sigma,\pi_4\sigma\}$ for
any permutation $\sigma\in \S_6$. In particular if
$\pi_0=(1,2,3,4,5,6)$ is the identity permutation, the taking
$\sigma=\pi_i^{-1}$ for $i=1,2,3,4$ will give other disallowed sets
permutations where one element of the set is the identity that are
essentially the same as $S$. Generally, this procedure gives 5 disallowed
sets from each one, which will result in a reduction in the number
by a factor of 5.  However, there is the possiblity that the set $S$
will actually be a group of permutations of order 5, in which case
$S\pi_i^{-1}=S$ for each $\pi_i\in S$.  In the present case, there
are four such groups of order 5, leaving a total of 290 sets which
are not groups, and thus there are 58 essential remaining cases.

Just as in the case of five observers of five events, in the case of
five observers of six events, where 294 of the 11,042,674,501 were
found to be disallowed, it is not known whether all of the remaining
11,042,674,207 cases can all occur. All that is known is that there
are allowable sign patterns of the $\alpha_j$'s (which are the
coefficients in the linear relation among the velocities), that might
make the pattern of permutations possible. 

Given the complexity associated with contructing point sets of
size five that realize sets of five permutations from $\S_5$, it seems
a bit daunting to try larger constructions.  However, in principle it is
clear how to generalize the prior construction. Starting with four
observers, one in the rest frame and the other three moving along
orthogonal axes, a set of six points in $3+1$-dimensional spacetime
can be constructed that realize four of the four permutations.  The
final observer frame can then be constucted if a set of linear
inequalities can be satisfied. The method used above now leads to two
``gap equations'', instead of one, namely
$g_4=\alpha_1\,g_1+\alpha_2\,g_2+\alpha_3\,g_3+\beta$ and
$g_5=\alpha^\prime_1\,g_1+\alpha^\prime_2\,g_2+\alpha^\prime_3\,g_3+\beta^\prime$
that need to be solved with $g_i>0$ for $i=1,2,3,4,5$.

Another phenomenon is that time reversal of a disallowed set should
yield another disallowed set. The time reversal of a set $S$ is
$\pi_rS\pi_r$ where $\pi_r=(6,5,4,3,2,1)$. The right multiplication
by $\pi_r$ accomplishes the time reversal itself, by reversing the
order of the observations for each observer.  The left multiplication
by $\pi_r$ simply relabels the observations. Of the 58 cases to consider,
8 were time reversal invariant, while 50 were not.  

The 8 time reversal invariant cases are listed below.  It is
interesting to note that one of these cases is also invariant under
taking the inverse of each permutation in the set.

\vskip 10pt
\settabs 20 \columns
\+  1:  &$\{(1,2,3,4,5,6),(1,4,6,5,3,2),(2,6,1,5,4,3),(4,3,2,6,1,5),(5,4,2,1,3,6)\}$&&&&&&&&&&&&&&&&&  \cr
\+  2:  &$\{(1,2,3,4,5,6),(1,4,6,5,3,2),(3,4,5,6,1,2),(5,4,2,1,3,6),(5,6,1,2,3,4)\}$&&&&&&&&&&&&&&&&&  \cr
\+  3:  &$\{(1,2,3,4,5,6),(1,5,6,4,2,3),(2,5,4,6,1,3),(4,5,3,1,2,6),(4,6,1,3,2,5)\}$&&&&&&&&&&&&&&&&&  \cr
\+  4:  &$\{(1,2,3,4,5,6),(1,5,6,4,2,3),(3,6,1,4,5,2),(4,5,3,1,2,6),(5,2,3,6,1,4)\}$&&&&&&&&&&&&&&&&& &&&self inverse  \cr
\+  5:  &$\{(1,2,3,4,5,6),(2,4,5,6,3,1),(3,4,5,6,1,2),(5,6,1,2,3,4),(6,4,1,2,3,5)\}$&&&&&&&&&&&&&&&&&  \cr
\+  6:  &$\{(1,2,3,4,5,6),(2,4,5,6,3,1),(4,3,5,1,6,2),(5,1,6,2,4,3),(6,4,1,2,3,5)\}$&&&&&&&&&&&&&&&&&  \cr
\+  7:  &$\{(1,2,3,4,5,6),(3,6,1,4,5,2),(4,2,3,6,5,1),(5,2,3,6,1,4),(6,2,1,4,5,3)\}$&&&&&&&&&&&&&&&&& \cr
\+  8:  &$\{(1,2,3,4,5,6),(4,1,6,3,5,2),(4,2,3,6,5,1),(5,2,4,1,6,3),(6,2,1,4,5,3)\}$&&&&&&&&&&&&&&&&& \cr

\vskip 20pt

None of the four groups of order five were time reversal
invariant. Thus, excluding the groups, there are 25 essentially
different cases that are not time reversal invariant . There are
essentially now only two different groups of order five.  None of
these cases are invariant under inversion. The complete results are given
in appendix 1.

\subhead
Five Observers of Five Events in 3-Dimensional Space
\endsubhead

Let $Q=\{\pi_0,\pi_1,\pi_2,\pi_3,\pi_4\}\subset \S_5$ be a set of
distinct event orderings of a set $\{E_1,E_2,E_3,E_4,E_5\}$ of five
(spacelike separated) events for 5 observers in 5 different inertial
frames. Without loss of generality, assume that $\pi_0$ is the
identity permutation in the rest frame, i.e. $\pi_0=(1,2,3,4,5)$.  In
this rest frame, event $E_i$ occurs at time $t_i$ and at spatial
coordinate $\tilde w_i$ for $i\in\{1,2,3,4,5\}$. Since the ordering of
events in the rest frame is the identity permutation, this is
equivalent to asserting that $t_1<t_2<t_3<t_4<t_5$. If observer $j$ is
moving at velocity $\tilde v^{(j)}$ from the perspective of the rest
frame (so $\tilde v^{(0)} = \bold 0$), then the time of event $E_i$ in 
frame $j$ is
$$T^{(j)}_{i}=\gamma^{(j)}\,\bigl(t_{i}-\tilde v^{(j)}\cdot\tilde w_{i}\bigr)$$
where $\gamma^{(j)}=\bigl(1-\tilde v\cdot\tilde v\bigr)^{-1/2}>0$.

In what follows it will be useful to deal with the ``relativized time'' 
of event $i$ for observer $j$ defined as
$$t^{(j)}_{i}=t_{i}-\tilde v^{(j)}\cdot\tilde w_{i}$$ 
which is the same as $T^{(j)}_{i}$ except for the Lorentz dilation
factor $\gamma^{(j)}$.  For the purposes of analyzing relative orders
of events being observed in frame $j$, the constant factor
$\gamma^{(j)}$ is irrelevant, and only makes the formulas more
complicated.  Also, in what follows, it will be mathematically
possible to have frame velocities larger than 1, which of course is
physically impossible.  It is always possible to divide all the
velocities by a large constant $N$ and then multiply the spatial
coordinates of all events by the same constant.  This way, in the
end, all the velocities can be made less than 1 with the observed
(relativized) event times being kept the same, and hence their ordering
in frame $j$ also is unchanged.

In addition to assuming that $\tilde v^{(0)} = \bold 0$, if three of the
remaining four velocity vectors span all of $\R^3$ (i.e. they don't
all come from a some smaller dimensional subspace) it is also possible
to change these three non-zero velocities to be any set of three
linearly independent velocities by making corresponding adjustments to
the coordinates of the events in spacetime.  Assuming $\tilde v^{(1)}$,
$\tilde v^{(2)}$, and $\tilde v^{(3)}$ are linearly independent,
if $\bigl\{\tilde V^{(1)},\tilde V^{(2)},\tilde V^{(3)}\bigr\}$ is any
other set of linearly independent velocity vectors, then there is an
invertible matrix $R$ such that $\tilde V^{(j)}=\tilde V^{(j)}\,R$ for
$j=1,2,3$. Now taking $W_{i} = R^{-1}\,w_{i}$ for all $i$ and 
$\tilde V^{(j)}=\tilde V^{(j)}\,R$ for all $j$ gives
$t_{i}-\tilde V^{(j)}\cdot\tilde W_{i} 
= t_{i}-\tilde v^{(j)}\cdot\tilde w_{i} = t^{(j)}_{i}$
for all $i$ and $j$. Thus without any loss of generality, it can be assumed
that $\tilde v^{(2)}$, and $\tilde v^{(3)}$ are orthogonal unit vectors in $\R^3$.

\subhead
Reducing the search to satisfying a single inequality
\endsubhead

What is desired is that 
$$t^{(j)}_{\pi_j(i+1)}>t^{(j)}_{\pi_j(i)}$$
for $i=1,2,3,4$. Now for $j=0,1,2,3$ let 
$$F_j:\{1,2,3,4,5\}\rightarrow\R$$
be a set of increasing functions.  It may be useful to write
$$h^{(j)}_i=F_j(i+1)-F_j(i)$$ 
for $j=0,1,2,3$ and $i=1,2,3,4$ with
the condition that the $h^{(j)}_i$'s are all positive, so that
$F_j(i) = F_j(1) + \sum_{l=1}^{i-1} h^{(j)}_l$
for $2\le i \le 5$. This in turn implies that
$$F_j(i_1) - F_j(i_0) 
= \sum_{l=i_0}^{i_1-1}h^{(j)}_l$$
for $1 \le i_0 < i_1 \le 5$.

Suppose that the spacetime coordinates of event $E_i$ in the rest frame are 
$$t_i=F_0(i)$$ 
and
$$\tilde w_i = (w_{i,1},w_{i,2},w_{i,2})
=\bigl(F_0(i)-F_1(\pi_1^{-1}(i)),
       F_0(i)-F_2(\pi_2^{-1}(i)),
       F_0(i)-F_3(\pi_3^{-1}(i))\bigr).$$
Then clearly $t^{(0)}_{i}=F_0(i)$ so 
$F_0(i+1)=t^{(0)}_{\pi_0(i+1)}>t^{(0)}_{\pi_0(i)}=F_0(i)$ for $i=1,2,3,4$ 
since $F_0$ is an increasing function. Now set
$$\eqalign{
\tilde v^{(1)} &= (1,0,0) \cr
\tilde v^{(2)} &= (0,1,0) \cr
\tilde v^{(3)} &= (0,0,1) \cr
}$$
so that 
$$t^{(1)}_{i}
=t_{i}-\tilde v^{(1)}\cdot\tilde w_{i}
=F_1(\pi_1^{-1}(i))
$$
for $i=1,2,3,4,5$, and therefore
$t^{(1)}_{\pi_1(i)} = F_1(i)$. Since $F_1$ is an increasing function,
it now follows that $t^{(1)}_{\pi_1(i+1)}>t^{(1)}_{\pi_1(i)}$
for $i=1,2,3,4$.
Similarly $t^{(2)}_{\pi_2(i)}=F_2(i)$ and
$t^{(3)}_{\pi_3(i)}=F_3(i)$ and therefore
$t^{(2)}_{\pi_2(i+1)}>t^{(2)}_{\pi_2(i)}$ and
$t^{(3)}_{\pi_3(i+1)}>t^{(3)}_{\pi_3(i)}$ for $i=1,2,3,4$ since $F_2$
and $F_3$ are also increasing functions. Thus, in general
$$t^{(j)}_{\pi_j(i+1)}-t^{(j)}_{\pi_j(i)}=F_j(i+1)-F_j(i)=h^{(j)}_i>0$$ 
for $i=1,2,3,4$ and $j=0,1,2,3$, by design

What remains is to consider the event ordering for the fifth observer,
whose velocity vector with respect to the rest frame is $\tilde v^{(4)}
=\tilde u=(u_1,u_2,u_3)$ with components that are to be determined. To this
end
$$\eqalign{t^{(4)}_{i}
&=t_{i}-\tilde v^{(4)}\cdot\tilde w_{i} \cr
&=F_0(i)-u_1\,(F_0(i)-F_1(\pi_1^{-1}(i)))
        -u_2\,(F_0(i)-F_2(\pi_2^{-1}(i))) \cr&\qquad
        -u_3\,(F_0(i)-F_3(\pi_3^{-1}(i))) \cr
&=F_0(i)\,(1-u_1-u_2-u_3)
     +u_1\,F_1(\pi_1^{-1}(i))
     +u_2\,F_2(\pi_2^{-1}(i))
     +u_3\,F_3(\pi_3^{-1}(i)) \cr
}$$
and therefore
$$\eqalign{t^{(4)}_{\pi_4(i)}
&=F_0(\pi_4(i))\,(1-u_1-u_2-u_3)
  +u_1\,F_1(\pi_1^{-1}(\pi_4(i))) \cr&\qquad
  +u_2\,F_2(\pi_2^{-1}(\pi_4(i)))
  +u_3\,F_3(\pi_3^{-1}(\pi_4(i))) \cr
}$$
for $i=1,2,3,4,5$, and note that the $t^{(4)}$'s are all linear functions
of the $u$'s. Now consider the gap between the values of the $t^{(4)}$'s,
by setting
$$\eqalign{
g_i
&=t^{(4)}_{\pi_4(i+1)}-t^{(4)}_{\pi_4(i)} \cr
&=(F_0(\pi_4(i+1))-F_0(\pi_4(i)))\bigl(1-u_1-u_2-u_3\bigr) 
\cr &\quad
     +u_1\,\bigl(F_1(\pi_1^{-1}(\pi_4(i+1))) - F_1(\pi_1^{-1}(\pi_4(i)))\bigr)
\cr &\quad
     +u_2\,\bigl(F_2(\pi_2^{-1}(\pi_4(i+1))) - F_2(\pi_2^{-1}(\pi_4(i)))\bigr)
\cr &\quad
     +u_3\,\bigl(F_3(\pi_3^{-1}(\pi_4(i+1))) - F_3(\pi_3^{-1}(\pi_4(i)))\bigr)
\cr}$$ 
for $i=1,2,3,4$. The conditions
$t^{(4)}_{\pi_4(i+1)}>t^{(4)}_{\pi_4(i)}$ for $i=1,2,3,4$ will be
satisfied if and only if all the $g_i$'s are positive.  However, the
$g_i$'s are linear functions of $u_1$, $u_2$, and $u_3$, so there is a
linear relation between the $g_i$'s. Assuming nonsingularity of the
linear system (which will be generically true, and can be guaranteed
by slight variations in the parameters), there is a linear
relationship
$$g_4=\alpha_1\,g_1+\alpha_2\,g_2+\alpha_3\,g_3+\beta$$
for some constants, $\alpha_1$, $\alpha_2$, $\alpha_3$, and $\beta$
that depend only on $F_0$, $F_1$, $F_2$, and $F_3$ and on $\pi_1$, $\pi_2$,
$\pi_3$, and $\pi_4$. If either $\beta$ is positive or the $\alpha$'s
are not all negative, then it is possible to take all the $g_i$'s to
be positive. Since the $g$'s are all linear functions of the $u$'s, it
will then be possible to find the desired velocity vector $v^{(4)}$ so
that the fifth observer will see the events in order $\pi_4$.

Write
$$g_i = b_i + \tilde a_i \cdot \tilde u$$
where $\tilde v^{(4)}=\tilde u=(u_1,u_2,u_3)$
and $\tilde a_i = (a_{i,1},a_{i,2},a_{i,3})$ with
$$\eqalign{a_{i,j} 
&= \bigl(F_j(\pi_j^{-1}(\pi_4(i+1))) - F_j(\pi_j^{-1}(\pi_4(i)))\bigr) 
  + (F_0(\pi_4(i)) - F_0(\pi_4(i+1))) \cr
&= w_{\pi_4(i+1),j} - w_{\pi_4(i),j} \cr
}$$
for $i=1,2,3,4$ and $j=1,2,3$ and
$$\eqalign{
b_i 
&= F_0(\pi_4(i+1)) - F_0(\pi_4(i)) \cr
&=  t_{\pi_4(i+1)} - t_{\pi_4(i)}  \cr
}$$
for $i=1,2,3,4$. Let $A$ be the $4\times3$ matrix with entries $(a_{i,j})$ as above,
and set $\tilde g = (g_1,g_2,g_3,g_4)$ and $\tilde b = (b_1,b_2,b_3,b_4)$
as row vectors. Then the above equations can be written as
$$\tilde g^{\tau} = A\,\tilde u^{\tau} + \tilde b^\tau.$$
Since $A$ has more rows than columns, it should generally
be possible to write the last row of $A$ as a linear combination of
the top three rows.  Write
$$
A=
\pmatrix
a_{1,1}&a_{1,2}&a_{1,3}\cr
a_{2,1}&a_{2,2}&a_{2,3}\cr
a_{3,1}&a_{3,2}&a_{3,3}\cr
a_{4,1}&a_{4,2}&a_{4,3}\cr
\endpmatrix
\qquad\hbox{and}\qquad
A_0=
\pmatrix
a_{1,1}&a_{2,1}&a_{3,1}\cr
a_{1,2}&a_{2,2}&a_{3,2}\cr
a_{1,3}&a_{2,3}&a_{3,3}\cr
\endpmatrix\
$$
and ask for a solution $\tilde \alpha = (\alpha_1\,\,\alpha_2\,\,\alpha_3)$ to
$$A_0\,\tilde \alpha^{\tau}=
\pmatrix
a_{1,1}&a_{2,1}&a_{3,1}\cr
a_{1,2}&a_{2,2}&a_{3,2}\cr
a_{1,3}&a_{2,3}&a_{3,3}\cr
\endpmatrix\,
\pmatrix
\alpha_1 \cr \alpha_2 \cr \alpha_3
\endpmatrix
=
\pmatrix
a_{4,1} \cr a_{4,2} \cr a_{4,3}
\endpmatrix
$$
so that
$\tilde \alpha^\prime
\,A=0$
where $\tilde \alpha^\prime = (\alpha_1\,\,\alpha_2\,\,\alpha_3\,-1)$. Then
$$\alpha_1\,g_1+\alpha_2\,g_2+\alpha_3\,g_3-g_4
=\tilde \alpha
\,\tilde g^{\tau} 
=\tilde \alpha^\prime
\,A\,\tilde u^{\tau}
  + \tilde \alpha^\prime
\,\tilde b^\tau
=\alpha_1\,b_1+\alpha_2\,b_2+\alpha_3\,b_3-b_4
$$
i.e.
$$g_4=\alpha_1\,g_1+\alpha_2\,g_2+\alpha_3\,g_3
       +(b_4-\alpha_1\,b_1-\alpha_2\,b_2-\alpha_3\,b_3)
$$
as the desired relation between the $g_i$'s. In particular
$\beta=b_4-\alpha_1\,b_1-\alpha_2\,b_2-\alpha_3\,b_3$.
Note that this depends on the nonsingularity of the system of equations,
but in general, by varying the $F_j$'s, this can be accomplished.

Of course, if the system of equations is singular, or if $\alpha_1$,
$\alpha_2$, $\alpha_3$, and $\beta$ are all negative, this program
will fail. In that case, a few options are available. Note that 
this whole approach started with picking out three of the four
elements of $\{\pi_1,\pi_2,\pi_3,\pi_4\}$ and setting the $x$, $y$,
and $z$ coordinates (in the rest frame) of the events 
$\{E_1,E_2,E_3,E_4,E_5\}$ to make those three permutations realizable.
Picking a different set of three permutations out of the four will
lead to a different set of equations to solve for the fourth velocity
vector, and may lead to a solution for the fourth permutation.

Another variation is to change the labelling on the events and pick
a different observer to be designated as the rest frame. This means
that the set of permutations $Q=\{\pi_0,\pi_1,\pi_2,\pi_3,\pi_4\}$ is
equivalent to the set $Q\sigma
=\{\pi_0\sigma,\pi_1\sigma,\pi_2\sigma,\pi_3\sigma,\pi_4\sigma\}$ for
any $\sigma$.  In particular, taking $\sigma=\pi_i^{-1}$ for
$i=1,2,3,4$ will give equivalent sets to $Q$ where one of the
permutations is $\pi_0$, and this will allow the entire theory above
to be applied. 

Another possible variation to note is that if a set $Q$ of
permutations is unachievable, then so is its time reversal $\pi^{(0)}
Q\pi^{(0)}$, where $\pi^{(0)}=(5,4,3,2,1)$. This may also lead to
additional disqualified sets.

Finally, there is always the option of varying the $F_j$'s as long as
they are all kept as increasing functions. This is equivalent to not
requiring the intervals between successive events in frames 1, 2, and 3
to be equally spaced. Note that the $a_{i,j}$'s themselves are simply
linear functions of the $h^{(j)}_i$'s and that the $h^{(j)}_i$'s are
all required to be positive. In particular the initial values $F_j(1)$
don't really matter since the $F_j$'s only appear as differences.

\subhead
Results of a computer search
\endsubhead

The ideas above were programmed and tested to see if every set of five
distinct permutations from $S_5$ was realizable. This is a
non-trivial, but also quite doable, search problem. This sort of
search problem is well within the range of computer technology, the
search space being of size 7,940,751.  In principle this can be cut down
to a search of size 1,588,155, although that was not actually done.

For each of the 7,940,751 choices of five distinct elements from $\S_5$ and for
each quadruple $(F_0,F_1,F_2,F_3)$ of increasing functions from
$\{1,2,3,4,5\}$ to $\R$, the above procedural description yields up to
twenty different linear relations of the form
$g_4=\alpha_1\,g_1+\alpha_2\,g_2+\alpha_3\,g_3+\beta$, and any one of
them that allows a solution with all the $g_i$'s positive will yield a
set of five points in Minkowski space and a set of five frame
velocities that will realize the desired set of permutations.  
The twenty possibilities come from five equivalent sets to the
set of permutations times four choices of three out of four
permutations to choose for orthogonal velocity vectors.

Note that one of the things that can go wrong here is that the set of
equations for the $\alpha_i$'s, which is of the form $\tilde a_4
=\alpha_1\,\tilde a_1 + \alpha_2\,\tilde a_2 + \alpha_3\,\tilde a_3$,
may in fact be singular, and therefore not solvable, which is why
there may actually be less than twenty different linear relations to
actually test.

The simplest choices of increasing functions are simply linear, i.e
$$F_j(x)=m_j\,x$$ 
where $m_j>0$ can be any positive real number.
There seems to be no systematic way of choosing a good set of $m_i$'s,
but if one set of $m_i$'s fails to yield any solvable linear relations,
another set of $m_i$'s can be readily tried, and that new set may well
succeed.  One somewhat bad choice, however, is taking $m_1=m_2=m_3=1$
since that doesn't really lead to twenty distinct linear relations due
to the symmetry of the problem, and in fact this leads only to 5
different linear relations.  

A run starting with $(m_0,m_1,m_2,m_3)=(1,300,200,300)$ yielded a
collection of only 333 permutation sets that were not realized.
\footnote{Some further analysis of the answer yielded an interesting
fact. For 160 of these cases, the linear system of equations for the
$\alpha_i$'s was singular regardless of how the permutations were
reordered or whether a equivalent set was chosen.  Also, these 160
cases were the same, regardless of how the $m_i$'s were changed ,
which clearly suggest that these 160 cases are intrinsically singular
and that there is an underlying identity.  For the remaining 173
cases, the system of equations for the $\alpha_i$'s were solvable, and
they were solvable regardless of how the permutations were ordered or
whether a equivalent set was used.  They just led to a set of equations
of the type $g_4=\alpha_1\,g_1+\alpha_2\,g_2+\alpha_3\,g_3+\beta$
where $\alpha_i<0$ for $i=1,2,3$ and $\beta<0$, and the signs of the
$\alpha$'s and $\beta$ didn't really depend on the choice of the
$m_i$'s.  However, of the 160 singular cases, only 60 appeared with
their time reverse, meaning that for 100 of the cases, where the
system was found to be singular, the time reverse was found to be
unachievable (and in particular also gave rise to a nonsingular
system), which means that this set of permutations was also
unachievable.}
Some experimentation with different $F_j$'s now helps a lot.  The following 
tabulated function
\vskip 10pt
\settabs 5 \columns
\+ $F_0(1) =    1.0$ & $F_0(2) =    2.0$ & $F_0(3) =    3.0$ & $F_0(4) =    4.0$ & $F_0(5) =    5.0$ &  \cr
\+ $F_1(1) =    1.0$ & $F_1(2) =    2.0$ & $F_1(3) =    3.0$ & $F_1(4) =    7.0$ & $F_1(5) =    8.0$ &  \cr
\+ $F_2(1) =    1.0$ & $F_2(2) =    2.0$ & $F_2(3) =   24.0$ & $F_2(4) =   25.0$ & $F_2(5) =   26.0$ &  \cr
\+ $F_3(1) =    1.0$ & $F_3(2) =   64.0$ & $F_3(3) =   65.0$ & $F_3(4) =   66.0$ & $F_3(5) =   67.0$ &  \cr
\vskip 10pt\noindent
realizes all but 3 cases of the 333 remaining permutation sets. A closer examination 
shows that each of these last three cases are cyclic subgroups of $\S_5$ of order 5,
and that two of these groups of order 5 are time reversal conjugates of each other.
Some further experimentation leads to the following tabulated function
\vskip 10pt
\settabs 5 \columns
\+ $F_0(1) =    1.0$ & $F_0(2) =    2.0$ & $F_0(3) =    3.0$ & $F_0(4) =    4.0$ & $F_0(5) =    5.0$ &  \cr
\+ $F_1(1) =    1.0$ & $F_1(2) =   23.0$ & $F_1(3) =   24.0$ & $F_1(4) =   25.0$ & $F_1(5) =   26.0$ &  \cr
\+ $F_2(1) =    4.0$ & $F_2(2) =    5.0$ & $F_2(3) =    6.0$ & $F_2(4) =    7.0$ & $F_2(5) =    8.0$ &  \cr
\+ $F_3(1) =    1.0$ & $F_3(2) =    2.0$ & $F_3(3) =   65.0$ & $F_3(4) =   66.0$ & $F_3(5) =   67.0$ &  \cr
\vskip 10pt\noindent
which realizes the two remaining permutation sets that are time reversals of 
each other.

This leaves one lone case still unresolved. This remaining case is the set
$$Q_0 = \{(1,2,3,4,5),(2,3,4,5,1),(3,4,5,1,2),(4,5,1,2,3),(5,1,2,3,4)\}$$
which is its own time reverse, and is actually invariant under inversion too,
since in fact this set of permutation is a cyclic group of order 5. 

So out of an initial set of 7,940,751 choices of five distinct
elements from $\S_5$ (containing the identity), all but at most one
set of permutations can be realized, i.e. for any such set of five
elements of $\S_5$, there are five events in three dimensional space
time and five observer frames such that each observer sees that the
five events in the specified permutation.

\subhead
Analysis of the final case
\endsubhead

It is useful to let $\pi=(2,3,4,5,1)$ be a generator of the cyclic
group under consideration.  Then the set of permutations to be
considered is just $Q_0 = \{\pi^0,\pi^1,\pi^2,\pi^3,\pi^4\}$, where
$\pi_j=\pi^j$.  If a set of representatives of $\Z/5\Z$ is taken to
be $\{1,2,3,4,5\}$ (rather than the more usual $\{0,1,2,3,4\}$), then
$\pi_j(i)=i+j\bmod5$. In particular $\pi_j^{-1}(\pi_4(i)) = i+4-j\bmod 5$
and $\pi_j^{-1}(\pi_4(i+1)) = i-j\bmod 5$, so
$$\eqalign{
F_j(\pi_j^{-1}(\pi_4(i))) - &F_j(\pi_j^{-1}(\pi_4(i+1))) 
= F_j(i+4-j\bmod5) - F_j(i-j\bmod5) \cr
&=\cases
h^{(j)}_1+h^{(j)}_2+h^{(j)}_3+h^{(j)}_4  & \hbox{if $i-j\equiv1\bmod5$;} \cr
-h^{(j)}_{(i-j\bmod 5)-1}                & \hbox{if $2\le (i-j\bmod 5) \le 5$.} \cr
\endcases
}$$
Similarly, $\pi_4(i) = i+4\bmod 5$ and $\pi_4(i+1) = i\bmod5$, so 
$$\eqalign{
F_0(\pi_4(i)) - &F_0(\pi_4(i+1)) 
= F_0(i+4\bmod5)-F_0(i\bmod5) \cr
&=\cases
h^{(0)}_1+h^{(0)}_2+h^{(0)}_3+h^{(0)}_4  & \hbox{if $i=1$;} \cr
-h^{(0)}_{i-1}                           & \hbox{if $2\le i \le 5$.} \cr
\endcases
}$$
Therefore
$$\eqalign{
&a_{i,j}
= \bigl(F_j(\pi_j^{-1}(\pi_4(i+1))) - F_j(\pi_j^{-1}(\pi_4(i)))\bigr) 
  + (F_0(\pi_4(i)) - F_0(\pi_4(i+1))) \cr
&=\cases
-h^{(j)}_1-h^{(j)}_2-h^{(j)}_3-h^{(j)}_4-h^{(0)}_{i-1} 
 & \hskip -.6pt \hbox{if $i-j\equiv1\bmod5$ and $2\le i \le 5$;} \cr
h^{(j)}_4 + h^{(0)}_1+h^{(0)}_2+h^{(0)}_3+h^{(0)}_4 
 & \hskip -.6pt \hbox{if $2\le (i-j\bmod 5) \le 5$ and $i=1$;} \cr
h^{(j)}_{(i-j\bmod 5)-1} -h^{(0)}_4 
 & \hskip -.6pt \hbox{if $2\le (i-j\bmod 5) \le 5$ and $2\le i \le 5$;} \cr
\endcases
}$$
where the case $i-j\equiv1\bmod5$ and $i=1$ has been omitted since 
this implies $j=0$, which is not relevant for what follows.
These are the entries of the matrix $4\times3$ $A$ which has entries $(a_{i,j})$ for 
$i=1,2,3,4$ and $j=1,2,3$, so that
$$\eqalign{A
&=\pmatrix
 {\textstyle h_4^{(1)}+h_1^{(0)}+h_2^{(0)}\atop\textstyle +h_3^{(0)}+h_4^{(0)} }
& {\textstyle h_3^{(2)}+h_1^{(0)}+h_2^{(0)}\atop\textstyle +h_3^{(0)}+h_4^{(0)} }
& {\textstyle h_2^{(3)}+h_1^{(0)}+h_2^{(0)}\atop\textstyle +h_3^{(0)}+h_4^{(0)} }\cr
 {\textstyle -h_1^{(1)}-h_2^{(1)}-h_3^{(1)}\atop\textstyle -h_4^{(1)}-h_1^{(0)} }
& h_4^{(2)}-h_1^{(0)} 
& h_3^{(3)}-h_1^{(0)} \cr
 h_1^{(1)}-h_2^{(0)} 
& {\textstyle -h_1^{(2)}-h_2^{(2)}-h_3^{(2)}\atop\textstyle -h_4^{(2)}-h_2^{(0)} }
& h_4^{(3)}-h_2^{(0)} \cr
 h_2^{(1)}-h_3^{(0)} 
& h_1^{(2)}-h_3^{(0)} 
& {\textstyle -h_1^{(3)}-h_2^{(3)}-h_3^{(3)}\atop\textstyle -h_4^{(3)}-h_3^{(0)} }\cr
\endpmatrix
\cr}$$
and 
$$\eqalign{
b_1 &= -h_1^{(0)}-h_2^{(0)}-h_3^{(0)}-h_4^{(0)} \cr
b_2 &= h_1^{(0)} \cr
b_3 &= h_2^{(0)} \cr
b_4 &= h_3^{(0)} \cr
}$$
since
$b_i = F_0(\pi_4(i+1)) - F_0(\pi_4(i))$
for $i=1,2,3,4$.

In general for a $4\times3$ matrix
$$
A=
\pmatrix
a_{1,1}&a_{1,2}&a_{1,3}\cr
a_{2,1}&a_{2,2}&a_{2,3}\cr
a_{3,1}&a_{3,2}&a_{3,3}\cr
a_{4,1}&a_{4,2}&a_{4,3}\cr
\endpmatrix
$$
there are four $3\times3$ matrices obtained by deleting a single row
$$\eqalign{
A_1=
\pmatrix
a_{2,1}&a_{2,2}&a_{2,3}\cr
a_{3,1}&a_{3,2}&a_{3,3}\cr
a_{4,1}&a_{4,2}&a_{4,3}\cr
\endpmatrix
&
A_2=
\pmatrix
a_{1,1}&a_{1,2}&a_{1,3}\cr
a_{3,1}&a_{3,2}&a_{3,3}\cr
a_{4,1}&a_{4,2}&a_{4,3}\cr
\endpmatrix
\cr
A_3=
\pmatrix
a_{1,1}&a_{1,2}&a_{1,3}\cr
a_{2,1}&a_{2,2}&a_{2,3}\cr
a_{4,1}&a_{4,2}&a_{4,3}\cr
\endpmatrix
&
A_4=
\pmatrix
a_{1,1}&a_{1,2}&a_{1,3}\cr
a_{2,1}&a_{2,2}&a_{2,3}\cr
a_{3,1}&a_{3,2}&a_{3,3}\cr
\endpmatrix
\cr}
$$
with the following general linear relation between the rows of $A$
$$\eqalign{0=
&\det(A_1)\cdot\bigl(a_{1,1},a_{1,2},a_{1,3}\bigr)
 -\det(A_2)\cdot\bigl(a_{2,1},a_{2,2},a_{2,3}\bigr) \cr&\quad
 +\det(A_3)\cdot\bigl(a_{3,1},a_{3,2},a_{3,3}\bigr)
 -\det(A_4)\cdot\bigl(a_{4,1},a_{4,2},a_{4,3}\bigr)
.\cr}$$
Writing $D_i = \det(A_i)$ for $i=1,2,3,4$, and setting
$$\eqalign{
\tilde a_1 &= \bigl(a_{1,1},a_{1,2},a_{1,3}\bigr) \cr
\tilde a_2 &= \bigl(a_{2,1},a_{2,2},a_{2,3}\bigr) \cr
\tilde a_3 &= \bigl(a_{3,1},a_{3,2},a_{3,3}\bigr) \cr
\tilde a_4 &= \bigl(a_{4,1},a_{4,2},a_{4,3}\bigr) \cr
}$$
so that if $D_4 \ne 0$ then
$$\tilde a_4 = \alpha_1\,\tilde a_1 + \alpha_2\,\tilde a_2 + \alpha_3\,\tilde a_3$$
with
$$\alpha_1 = D_1/D_4 \qquad \alpha_2 = -D_2/D_4 \qquad \alpha_3=D_3/D_4.$$

It is a straightforward (but uncomfortably large) algebraic calculation to find the
$D_k$'s as functions of the $h_i^{(j)}$'s, and the details of this have been relegated
to appendix 2. However, the point of the whole computation is to show that  
$D_1,D_3<0$ and $D_2,D_4>0$. Therefore $\alpha_1,\alpha_2,\alpha_3<0$. Also
$$\beta = b_4 - \alpha_1\,b_1 - \alpha_2\,b_2 - \alpha_3\,b_3$$ 
and the computation of $\beta$ (also relegated to appendix 2), shows that $\beta < 0$, as well.

Since $\alpha_i<0$ for $i=1,2,3$ and $\beta < 0$, it follows that
$g_4=\alpha_1\,g_1+\alpha_2\,g_2+\alpha_3\,g_3+\beta < 0$ if $g_i>0$
for $i=1,2,3$, and therefore there is no velocity vector $\tilde
v^{(4)} =\tilde u=(u_1,u_2,u_3)$ that realizes the fourth permutation.

\proclaim{Theorem}
Let $Q\subset\S_5$ with $\card{Q}=5$ and $\pi_0\in Q$. If $Q\ne Q_0$, then there exist a
set of 5 points in spacetime and a set of 5 inertial reference frames
that realize the set $Q$. If $Q=Q_0$, then such a set of 5 points in spacetime and a set of
5 reference frames does not exist.
\endproclaim

\subhead
Implications for five reference frames and six points
\endsubhead

With a single forbidden configuation of five permutation on five
elements, it is easy to construct sets of five forbidden permutations
on six elements by simply slipping in the sixth element anywhere in
each of the five permutations comprising $Q_0$.

In principle this new event is viewed as an insertion of a new number
in each extant permutation giving now a set of five elements of
$\S_6$. This insertion can occur in any one of six possible points for
each $\pi^i$, and all possibilities must be considered. There are
$6^5$ such possibilities, however for each possibility, it is
necessary to convert it to "standard" form where one of the
permutations is the identity.  There five ways to do this
corresponding to five different choices of which observer sees the
identity permutation, making a total of $5\cdot6^5=38,880$ possible
cases to consider.  Note that when this is done, it may be that some
of cases are the same.  The way that duplicates are eliminated is by
sorting the set of cases.

When this is done, and all the duplication is elimated, there are 7676
subsets of $\S_6$ of cardinality 5 that contain the identity
permutation remaining.  All of these represent forbidden
configurations in addition to the 294 forbidden configurations already
found by considering sign patterns.  It is interesting that this new
collection of 7676 forbidden sets is completely disjoint from the
prior collection of 294 forbidden sets.

As discussed above, if $S =
\{\pi_0,\pi_1,\pi_2,\pi_3,\pi_4\}$ is any disallowed pattern, then so
is $S\sigma$ and taking
$\sigma=\pi_i^{-1}$ for $i=1,2,3,4$ will give other disallowed sets
permutations where one element of the set is the identity.  These
the same as $S$, and generally, this procedure gives 5 disallowed
sets from each one, however there is the possiblity that the set $S$
will actually be a group of permutations of order 5 contained in $S_6$.
In the present case, of the 7676 forbidden cases coming from $Q_0$, there
are six such groups of order 5, leaving a total of 7670 sets which are
not groups, and a total of 1534 essential remaining cases. Of the 1540
essentially different cases, only six were invariant under time reversal
(and interestingly one of these is actually invariant under permutation
inversion), while 1534 were not. The total number of cases to consider
then is $6+1534/2=773$.  It would be too long, and probably not of great
interest to list all 773 cases, since they are all easily constructed
from the basic $Q_0$ set of five elements of $\S_5$.  It may be of some
interest to list the six time reversal invariant sets.  They are:
\vskip 9pt
\settabs 24 \columns
\+  1:  &$\{(1,2,3,4,5,6),(3,4,5,1,6,2),(4,5,1,6,2,3),(5,1,6,2,3,4),(6,2,3,4,5,1)\}$&&&&&&&&&&&&&&&&&  \cr
\+  2:  &$\{(1,2,3,4,5,6),(3,4,5,1,6,2),(4,5,6,1,2,3),(5,1,6,2,3,4),(6,2,3,4,5,1)\}$&&&&&&&&&&&&&&&&&  \cr
\+  3:  &$\{(1,2,3,4,5,6),(2,3,4,5,6,1),(3,4,5,6,1,2),(5,6,1,2,3,4),(6,1,2,3,4,5)\}$&&&&&&&&&&&&&&&&&&&&&& self inverse  \cr
\+  4:  &$\{(1,2,3,4,5,6),(2,4,3,5,6,1),(4,3,5,6,1,2),(5,6,1,2,4,3),(6,1,2,4,3,5)\}$&&&&&&&&&&&&&&&&& \cr
\+  5:  &$\{(1,2,3,4,5,6),(2,4,3,5,6,1),(3,4,5,6,1,2),(5,6,1,2,3,4),(6,1,2,4,3,5)\}$&&&&&&&&&&&&&&&&& \cr
\+  6:  &$\{(1,2,3,4,5,6),(2,3,4,5,6,1),(4,3,5,6,1,2),(5,6,1,2,4,3),(6,1,2,3,4,5)\}$&&&&&&&&&&&&&&&&& \cr
\vskip 9pt

A seemingly much more complicated question is whether there are other
unrealizable subsets of $\S_6$ of size five that are unrealizable
that are not included in the $7676+294=7970$ cases already considered.

\head
Conclusion
\endhead

This note represents an initial attempt to capture the combinatorial
nature of special relativity.  It is an interesting question whether
starting with the combinatorial restrictions imposed by special
relativity eventually lead to Minkowsky space-time.

In 3+1-dimensional spacetime with 5 observers of 5 spacelike separated
events, it is striking that of the 7,940,751 possible cases to
consider, all but exactly one are realizable. The proof is based on a
computer construction of the 7,940,750 realizable cases along with a
computer algebra computation to show that the remaining case is
impossible. One might hope that a more conceptual proof could be
found. With 5 observers of 6 spacelike separated events, along with
extending the impossible case of 5 events to an arbitrary sixth event,
a counting of sign restrictions shows there to be at least 294
additional unrealizable sets of permutations, also by a computer
search.
As a concluding comment, given the minimal non-realizable sets that
have been shown here for $1+1$ relativity and for $3+1$ relativity, it
might be reasonable to conjecture that in $n+1$ relativity, there is
exactly one non-realizable set of $n+1$ events for $n+1$ observers,
specifically the set based on the cyclic permutation.

\vfill\break
\subhead
Appendix 1: Five Observers of Six Events in 3-Dimensional Space: All Cases
\endsubhead

All 35 cases are listed below.  The first eight entries are the time
invariant sets and the remainder of the listing pairs up each set with
its time reversal.

\vskip 10pt
\settabs 20 \columns
\+  1:  &$\{(1,2,3,4,5,6),(1,4,6,5,3,2),(2,6,1,5,4,3),(4,3,2,6,1,5),(5,4,2,1,3,6)\}$&&&&&&&&&&&&&&&&&  \cr
\+  2:  &$\{(1,2,3,4,5,6),(1,4,6,5,3,2),(3,4,5,6,1,2),(5,4,2,1,3,6),(5,6,1,2,3,4)\}$&&&&&&&&&&&&&&&&&  \cr
\+  3:  &$\{(1,2,3,4,5,6),(1,5,6,4,2,3),(2,5,4,6,1,3),(4,5,3,1,2,6),(4,6,1,3,2,5)\}$&&&&&&&&&&&&&&&&&  \cr
\+  4:  &$\{(1,2,3,4,5,6),(1,5,6,4,2,3),(3,6,1,4,5,2),(4,5,3,1,2,6),(5,2,3,6,1,4)\}$&&&&&&&&&&&&&&&&& &&&self inverse  \cr
\+  5:  &$\{(1,2,3,4,5,6),(2,4,5,6,3,1),(3,4,5,6,1,2),(5,6,1,2,3,4),(6,4,1,2,3,5)\}$&&&&&&&&&&&&&&&&&  \cr
\+  6:  &$\{(1,2,3,4,5,6),(2,4,5,6,3,1),(4,3,5,1,6,2),(5,1,6,2,4,3),(6,4,1,2,3,5)\}$&&&&&&&&&&&&&&&&&  \cr
\+  7:  &$\{(1,2,3,4,5,6),(3,6,1,4,5,2),(4,2,3,6,5,1),(5,2,3,6,1,4),(6,2,1,4,5,3)\}$&&&&&&&&&&&&&&&&& \cr
\+  8:  &$\{(1,2,3,4,5,6),(4,1,6,3,5,2),(4,2,3,6,5,1),(5,2,4,1,6,3),(6,2,1,4,5,3)\}$&&&&&&&&&&&&&&&&& \cr

\vskip 20pt 
\+ 9a:
&$\{(1,2,3,4,5,6),(1,4,6,5,3,2),(2,4,6,5,1,3),(5,6,1,2,3,4),(6,3,2,1,4,5)\}$&&&&&&&&&&&&&&&&& &&&&&&\cr 
\+ 9b: &$\{(1,2,3,4,5,6),(2,3,6,5,4,1),(3,4,5,6,1,2),(4,6,2,1,3,5),(5,4,2,1,3,6)\}$&&&&&&&&&&&&&&&&& &&&&&&\cr 
\vskip 5pt 
\+ 10a: &$\{(1,2,3,4,5,6),(1,4,6,5,3,2),(3,4,5,6,1,2),(5,6,1,2,3,4),(6,3,2,1,4,5)\}$&&&&&&&&&&&&&&&&& &&&&&&\cr 
\+ 10b: &$\{(1,2,3,4,5,6),(2,3,6,5,4,1),(3,4,5,6,1,2),(5,4,2,1,3,6),(5,6,1,2,3,4)\}$&&&&&&&&&&&&&&&&& &&&&&&\cr 
\vskip 5pt 
\+ 11a: &$\{(1,2,3,4,5,6),(1,4,6,5,3,2),(3,5,4,1,6,2),(4,3,2,6,1,5),(6,3,1,2,5,4)\}$&&&&&&&&&&&&&&&&& &&&&&&\cr 
\+ 11b: &$\{(1,2,3,4,5,6),(2,6,1,5,4,3),(3,2,5,6,4,1),(5,1,6,3,2,4),(5,4,2,1,3,6)\}$&&&&&&&&&&&&&&&&& &&&&&&\cr 
\vskip 5pt 
\+ 12a: &$\{(1,2,3,4,5,6),(1,4,6,5,3,2),(3,6,1,5,2,4),(4,3,2,6,1,5),(5,3,4,1,2,6)\}$&&&&&&&&&&&&&&&&& &&&&&&\cr 
\+ 12b: &$\{(1,2,3,4,5,6),(1,5,6,3,4,2),(2,6,1,5,4,3),(3,5,2,6,1,4),(5,4,2,1,3,6)\}$&&&&&&&&&&&&&&&&& &&&&&&\cr 
\vskip 5pt 
\+ 13a: &$\{(1,2,3,4,5,6),(1,4,6,5,3,2),(3,6,4,1,2,5),(5,3,1,6,2,4),(6,2,1,5,3,4)\}$&&&&&&&&&&&&&&&&& &&&&&&\cr 
\+ 13b: &$\{(1,2,3,4,5,6),(2,5,6,3,1,4),(3,4,2,6,5,1),(3,5,1,6,4,2),(5,4,2,1,3,6)\}$&&&&&&&&&&&&&&&&& &&&&&&\cr 
\vskip 5pt 
\+ 14a: &$\{(1,2,3,4,5,6),(1,4,6,5,3,2),(4,2,5,6,1,3),(5,3,4,1,2,6),(5,6,1,2,3,4)\}$&&&&&&&&&&&&&&&&& &&&&&&\cr 
\+ 14b: &$\{(1,2,3,4,5,6),(1,5,6,3,4,2),(3,4,5,6,1,2),(4,6,1,2,5,3),(5,4,2,1,3,6)\}$&&&&&&&&&&&&&&&&& &&&&&&\cr 
\vskip 5pt 
\+ 15a: &$\{(1,2,3,4,5,6),(1,4,6,5,3,2),(4,2,6,3,1,5),(5,2,4,1,6,3),(6,2,1,5,3,4)\}$&&&&&&&&&&&&&&&&& &&&&&&\cr 
\+ 15b: &$\{(1,2,3,4,5,6),(2,6,4,1,5,3),(3,4,2,6,5,1),(4,1,6,3,5,2),(5,4,2,1,3,6)\}$&&&&&&&&&&&&&&&&& &&&&&&\cr 
\vskip 5pt 
\+ 16a: &$\{(1,2,3,4,5,6),(1,4,6,5,3,2),(5,2,4,1,6,3),(5,3,1,6,2,4),(6,2,1,5,3,4)\}$&&&&&&&&&&&&&&&&& &&&&&&\cr 
\+ 16b: &$\{(1,2,3,4,5,6),(3,4,2,6,5,1),(3,5,1,6,4,2),(4,1,6,3,5,2),(5,4,2,1,3,6)\}$&&&&&&&&&&&&&&&&& &&&&&&\cr 
\vskip 5pt 
\+ 17a: &$\{(1,2,3,4,5,6),(1,5,4,6,3,2),(2,5,6,1,4,3),(4,5,2,3,1,6),(4,6,2,1,3,5)\}$&&&&&&&&&&&&&&&&& &&&&&&\cr 
\+ 17b: &$\{(1,2,3,4,5,6),(1,6,4,5,2,3),(2,4,6,5,1,3),(4,3,6,1,2,5),(5,4,1,3,2,6)\}$&&&&&&&&&&&&&&&&& &&&&&&\cr 
\vskip 5pt 
\+ 18a: &$\{(1,2,3,4,5,6),(1,5,4,6,3,2),(2,5,6,1,4,3),(4,6,2,1,3,5),(6,3,1,2,5,4)\}$&&&&&&&&&&&&&&&&& &&&&&&\cr 
\+ 18b: &$\{(1,2,3,4,5,6),(2,4,6,5,1,3),(3,2,5,6,4,1),(4,3,6,1,2,5),(5,4,1,3,2,6)\}$&&&&&&&&&&&&&&&&& &&&&&&\cr 
\vskip 5pt 
\+ 19a: &$\{(1,2,3,4,5,6),(1,5,4,6,3,2),(3,4,5,6,1,2),(4,6,2,1,3,5),(6,3,1,2,5,4)\}$&&&&&&&&&&&&&&&&& &&&&&&\cr 
\+ 19b: &$\{(1,2,3,4,5,6),(2,4,6,5,1,3),(3,2,5,6,4,1),(5,4,1,3,2,6),(5,6,1,2,3,4)\}$&&&&&&&&&&&&&&&&& &&&&&&\cr 
\vskip 5pt 
\+ 20a: &$\{(1,2,3,4,5,6),(1,5,4,6,3,2),(3,4,5,6,1,2),(5,3,2,1,6,4),(6,3,1,2,5,4)\}$&&&&&&&&&&&&&&&&& &&&&&&\cr 
\+ 20b: &$\{(1,2,3,4,5,6),(3,1,6,5,4,2),(3,2,5,6,4,1),(5,4,1,3,2,6),(5,6,1,2,3,4)\}$&&&&&&&&&&&&&&&&& &&&&&&\cr 
\vskip 5pt 
\+ 21a: &$\{(1,2,3,4,5,6),(1,5,6,3,4,2),(2,5,6,1,4,3),(3,6,2,1,5,4),(5,3,2,4,1,6)\}$&&&&&&&&&&&&&&&&& &&&&&&\cr 
\+ 21b: &$\{(1,2,3,4,5,6),(1,6,3,5,4,2),(3,2,6,5,1,4),(4,3,6,1,2,5),(5,3,4,1,2,6)\}$&&&&&&&&&&&&&&&&& &&&&&&\cr 
\vskip 5pt 
\+ 22a: &$\{(1,2,3,4,5,6),(1,5,6,3,4,2),(2,5,6,1,4,3),(3,6,2,1,5,4),(6,4,1,2,3,5)\}$&&&&&&&&&&&&&&&&& &&&&&&\cr 
\+ 22b: &$\{(1,2,3,4,5,6),(2,4,5,6,3,1),(3,2,6,5,1,4),(4,3,6,1,2,5),(5,3,4,1,2,6)\}$&&&&&&&&&&&&&&&&& &&&&&&\cr 
\vskip 5pt 
\+ 23a: &$\{(1,2,3,4,5,6),(1,5,6,3,4,2),(2,5,6,1,4,3),(4,6,1,2,5,3),(6,3,2,1,4,5)\}$&&&&&&&&&&&&&&&&& &&&&&&\cr 
\+ 23b: &$\{(1,2,3,4,5,6),(2,3,6,5,4,1),(4,2,5,6,1,3),(4,3,6,1,2,5),(5,3,4,1,2,6)\}$&&&&&&&&&&&&&&&&& &&&&&&\cr 
\vskip 5pt 
\+ 24a: &$\{(1,2,3,4,5,6),(1,5,6,3,4,2),(3,4,5,6,1,2),(4,5,2,1,6,3),(6,4,1,2,3,5)\}$&&&&&&&&&&&&&&&&& &&&&&&\cr 
\+ 24b: &$\{(1,2,3,4,5,6),(2,4,5,6,3,1),(4,1,6,5,2,3),(5,3,4,1,2,6),(5,6,1,2,3,4)\}$&&&&&&&&&&&&&&&&& &&&&&&\cr 
\vskip 5pt 
\+ 25a: &$\{(1,2,3,4,5,6),(1,5,6,4,2,3),(2,5,4,6,1,3),(4,6,1,3,2,5),(6,2,3,1,5,4)\}$&&&&&&&&&&&&&&&&& &&&&&&\cr 
\+ 25b: &$\{(1,2,3,4,5,6),(2,5,4,6,1,3),(3,2,6,4,5,1),(4,5,3,1,2,6),(4,6,1,3,2,5)\}$&&&&&&&&&&&&&&&&& &&&&&&\cr 
\vskip 5pt 
\+ 26a: &$\{(1,2,3,4,5,6),(1,5,6,4,2,3),(2,5,4,6,1,3),(5,3,2,1,6,4),(6,2,3,1,5,4)\}$&&&&&&&&&&&&&&&&& &&&&&&\cr 
\+ 26b: &$\{(1,2,3,4,5,6),(3,1,6,5,4,2),(3,2,6,4,5,1),(4,5,3,1,2,6),(4,6,1,3,2,5)\}$&&&&&&&&&&&&&&&&& &&&&&&\cr 
\vskip 5pt 
\+ 27a: &$\{(1,2,3,4,5,6),(1,6,3,5,4,2),(2,6,1,5,4,3),(3,6,2,4,1,5),(5,3,2,1,6,4)\}$&&&&&&&&&&&&&&&&& &&&&&&\cr 
\+ 27b: &$\{(1,2,3,4,5,6),(2,6,3,5,1,4),(3,1,6,5,4,2),(4,3,2,6,1,5),(5,3,2,4,1,6)\}$&&&&&&&&&&&&&&&&& &&&&&&\cr 
\vskip 5pt 
\+ 28a: &$\{(1,2,3,4,5,6),(1,6,3,5,4,2),(2,6,5,1,3,4),(3,6,2,4,1,5),(4,6,1,2,5,3)\}$&&&&&&&&&&&&&&&&& &&&&&&\cr 
\+ 28b: &$\{(1,2,3,4,5,6),(2,6,3,5,1,4),(3,4,6,2,1,5),(4,2,5,6,1,3),(5,3,2,4,1,6)\}$&&&&&&&&&&&&&&&&& &&&&&&\cr 
\vskip 5pt 
\+ 29a: &$\{(1,2,3,4,5,6),(2,3,6,5,4,1),(3,1,6,5,4,2),(5,6,1,2,3,4),(6,4,1,2,3,5)\}$&&&&&&&&&&&&&&&&& &&&&&&\cr 
\+ 29b: &$\{(1,2,3,4,5,6),(2,4,5,6,3,1),(3,4,5,6,1,2),(5,3,2,1,6,4),(6,3,2,1,4,5)\}$&&&&&&&&&&&&&&&&& &&&&&&\cr 
\vskip 5pt 
\+ 30a: &$\{(1,2,3,4,5,6),(2,3,6,5,4,1),(3,4,5,6,1,2),(5,6,1,2,3,4),(6,4,1,2,3,5)\}$&&&&&&&&&&&&&&&&& &&&&&&\cr 
\+ 30b: &$\{(1,2,3,4,5,6),(2,4,5,6,3,1),(3,4,5,6,1,2),(5,6,1,2,3,4),(6,3,2,1,4,5)\}$&&&&&&&&&&&&&&&&& &&&&&&\cr 
\vskip 5pt 
\+ 31a: &$\{(1,2,3,4,5,6),(2,4,5,6,3,1),(3,6,2,4,1,5),(4,3,5,1,6,2),(5,3,2,1,6,4)\}$&&&&&&&&&&&&&&&&& &&&&&&\cr 
\+ 31b: &$\{(1,2,3,4,5,6),(2,6,3,5,1,4),(3,1,6,5,4,2),(5,1,6,2,4,3),(6,4,1,2,3,5)\}$&&&&&&&&&&&&&&&&& &&&&&&\cr 
\vskip 5pt 
\+ 32a: &$\{(1,2,3,4,5,6),(2,4,5,6,3,1),(4,6,1,3,2,5),(5,1,4,6,2,3),(6,2,1,5,3,4)\}$&&&&&&&&&&&&&&&&& &&&&&&\cr 
\+ 32b: &$\{(1,2,3,4,5,6),(2,5,4,6,1,3),(3,4,2,6,5,1),(4,5,1,3,6,2),(6,4,1,2,3,5)\}$&&&&&&&&&&&&&&&&& &&&&&&\cr 
\vskip 5pt 
\+ 33a: &$\{(1,2,3,4,5,6),(2,4,6,3,5,1),(3,6,2,1,5,4),(4,3,1,6,5,2),(6,1,4,2,5,3)\}$&&&&&&&&&&&&&&&&& &&&group of order 5 &&&&&&\cr 
\+ 33b: &$\{(1,2,3,4,5,6),(3,2,6,5,1,4),(4,2,5,3,6,1),(5,2,1,6,4,3),(6,2,4,1,3,5)\}$&&&&&&&&&&&&&&&&& &&&group of order 5 &&&&&&\cr 
\vskip 5pt 
\+ 34a: &$\{(1,2,3,4,5,6),(2,4,6,3,5,1),(3,6,4,1,2,5),(4,1,6,5,2,3),(6,2,1,5,3,4)\}$&&&&&&&&&&&&&&&&& &&&&&&\cr 
\+ 34b: &$\{(1,2,3,4,5,6),(2,5,6,3,1,4),(3,4,2,6,5,1),(4,5,2,1,6,3),(6,2,4,1,3,5)\}$&&&&&&&&&&&&&&&&& &&&&&&\cr 
\vskip 5pt 
\+ 35a: &$\{(1,2,3,4,5,6),(2,6,4,1,5,3),(3,4,2,6,5,1),(4,1,6,3,5,2),(6,3,1,2,5,4)\}$&&&&&&&&&&&&&&&&& &&&group of order 5 &&&&&&\cr 
\+ 35b: &$\{(1,2,3,4,5,6),(3,2,5,6,4,1),(4,2,6,3,1,5),(5,2,4,1,6,3),(6,2,1,5,3,4)\}$&&&&&&&&&&&&&&&&& &&&group of order 5 &&&&&&\cr 
\vskip 5pt

\vfill\break
\subhead
Appendix 2: Calculation of $\alpha_i$'s and $\beta$
\endsubhead

For a $4\times3$ matrix
$$
A=
\pmatrix
a_{1,1}&a_{1,2}&a_{1,3}\cr
a_{2,1}&a_{2,2}&a_{2,3}\cr
a_{3,1}&a_{3,2}&a_{3,3}\cr
a_{4,1}&a_{4,2}&a_{4,3}\cr
\endpmatrix
$$
consider the four $3\times3$ matrices obtained by deleting a single row
$$\eqalign{
A_1=
\pmatrix
a_{2,1}&a_{2,2}&a_{2,3}\cr
a_{3,1}&a_{3,2}&a_{3,3}\cr
a_{4,1}&a_{4,2}&a_{4,3}\cr
\endpmatrix
&\qquad
A_2=
\pmatrix
a_{1,1}&a_{1,2}&a_{1,3}\cr
a_{3,1}&a_{3,2}&a_{3,3}\cr
a_{4,1}&a_{4,2}&a_{4,3}\cr
\endpmatrix
\cr
A_3=
\pmatrix
a_{1,1}&a_{1,2}&a_{1,3}\cr
a_{2,1}&a_{2,2}&a_{2,3}\cr
a_{4,1}&a_{4,2}&a_{4,3}\cr
\endpmatrix
&\qquad
A_4=
\pmatrix
a_{1,1}&a_{1,2}&a_{1,3}\cr
a_{2,1}&a_{2,2}&a_{2,3}\cr
a_{3,1}&a_{3,2}&a_{3,3}\cr
\endpmatrix
\cr}
$$
There is the following general linear relation between the rows of $A$
$$\eqalign{0=
& \det(A_1)\cdot\bigl(a_{1,1},a_{1,2},a_{1,3}\bigr)
 -\det(A_2)\cdot\bigl(a_{2,1},a_{2,2},a_{2,3}\bigr) \cr&\quad
 +\det(A_3)\cdot\bigl(a_{3,1},a_{3,2},a_{3,3}\bigr)
 -\det(A_4)\cdot\bigl(a_{4,1},a_{4,2},a_{4,3}\bigr)
}$$
Writing $D_i = \det(A_i)$ for $i=1,2,3,4$, and setting
$$\eqalign{
\tilde a_1 &= \bigl(a_{1,1},a_{1,2},a_{1,3}\bigr) \cr
\tilde a_2 &= \bigl(a_{2,1},a_{2,2},a_{2,3}\bigr) \cr
\tilde a_3 &= \bigl(a_{3,1},a_{3,2},a_{3,3}\bigr) \cr
\tilde a_4 &= \bigl(a_{4,1},a_{4,2},a_{4,3}\bigr) \cr
}$$
so that if $D_4 \ne 0$ then
$$\tilde a_4 = \alpha_1\,\tilde a_1 + \alpha_2\,\tilde a_2 + \alpha_3\,\tilde a_3$$
with
$$\alpha_1 = D_1/D_4 \qquad \alpha_2 = -D_2/D_4 \qquad \alpha_3=D_3/D_4$$

For the set 
$$Q_0=\{(1,2,3,4,5),(2,3,4,5,1),(3,4,5,1,2),(4,5,1,2,3),(5,1,2,3,4)\}$$
the $4\times3$ matrix is 
$$\scriptstyle A=\pmatrix 
\scriptstyle h_4^{(1)}+h_1^{(0)}+h_2^{(0)}+h_3^{(0)}+h_4^{(0)} 
& \scriptstyle h_3^{(2)}+h_1^{(0)}+h_2^{(0)}+h_3^{(0)}+h_4^{(0)} 
& \scriptstyle h_2^{(3)}+h_1^{(0)}+h_2^{(0)}+h_3^{(0)}+h_4^{(0)} \cr
\scriptstyle  -h_1^{(1)}-h_2^{(1)}-h_3^{(1)}-h_4^{(1)}-h_1^{(0)} 
& \scriptstyle h_4^{(2)}-h_1^{(0)} 
& \scriptstyle h_3^{(3)}-h_1^{(0)} \cr
 \scriptstyle h_1^{(1)}-h_2^{(0)} 
& \scriptstyle -h_1^{(2)}-h_2^{(2)}-h_3^{(2)}-h_4^{(2)}-h_2^{(0)} 
& \scriptstyle h_4^{(3)}-h_2^{(0)} \cr
 \scriptstyle h_2^{(1)}-h_3^{(0)} 
& \scriptstyle h_1^{(2)}-h_3^{(0)} 
& \scriptstyle -h_1^{(3)}-h_2^{(3)}-h_3^{(3)}-h_4^{(3)}-h_3^{(0)} \cr
\endpmatrix$$
and 
$$\eqalign{
b_1 &= -h_1^{(0)}-h_2^{(0)}-h_3^{(0)}-h_4^{(0)} \cr
b_2 &= h_1^{(0)} \cr
b_3 &= h_2^{(0)} \cr
b_4 &= h_3^{(0)} \cr
}$$
Writing $D_i = \det(A_i)$ for $i=1,2,3,4$, the computations below show that
$D_1,D_3<0$ and $D_2,D_4>0$. Therefore $\alpha_1,\alpha_2,\alpha_3<0$. Also
$$\beta = b_4 - \alpha_1\,b_1 - \alpha_2\,b_2 - \alpha_3\,b_3$$ 
and a computation of $\beta$ is also given below, which shows that $\beta < 0$.
Putting all this together demonstrates that the permutation set $Q_0$ is not
realizable.

\break
The following results were obtained by the computer algebra system magma:

$$\scriptstyle  A_4^\tau= \pmatrix
\scriptstyle h^{(0)}_{1} + h^{(0)}_{2} + h^{(0)}_{3} + h^{(0)}_{4} + h^{(1)}_{4} &
\scriptstyle     h^{(0)}_{1} + h^{(0)}_{2} + h^{(0)}_{3} + h^{(2)}_{3} + h^{(0)}_{4} &
\scriptstyle     h^{(0)}_{1} + h^{(0)}_{2} + h^{(3)}_{2} + h^{(0)}_{3} + h^{(0)}_{4} \cr
\scriptstyle -h^{(0)}_{1} - h^{(1)}_{1} - h^{(1)}_{2} - h^{(1)}_{3} - h^{(1)}_{4} &
\scriptstyle    -h^{(0)}_{1} + h^{(2)}_{4}&
\scriptstyle    -h^{(0)}_{1} + h^{(3)}_{3} \cr
\scriptstyle h^{(1)}_{1} - h^{(0)}_{2}   &
\scriptstyle    -h^{(2)}_{1} - h^{(0)}_{2} - h^{(2)}_{2} - h^{(2)}_{3} - h^{(2)}_{4} &
\scriptstyle    -h^{(0)}_{2} + h^{(3)}_{4} \cr
\endpmatrix
$$

$$\eqalign{
\scriptstyle D_4&\scriptstyle = \det(A_4) \cr&= \scriptstyle
    h^{(0)}_{1}\,h^{(1)}_{1}\,h^{(2)}_{1} + h^{(0)}_{1}\,h^{(1)}_{1}\,h^{(2)}_{2} +
    h^{(0)}_{1}\,h^{(1)}_{1}\,h^{(3)}_{2} + h^{(0)}_{1}\,h^{(1)}_{1}\,h^{(3)}_{3} +
    h^{(0)}_{1}\,h^{(1)}_{1}\,h^{(3)}_{4} + h^{(0)}_{1}\,h^{(2)}_{1}\,h^{(1)}_{2} \cr&\,\,\scriptstyle+
    h^{(0)}_{1}\,h^{(2)}_{1}\,h^{(3)}_{2} + h^{(0)}_{1}\,h^{(2)}_{1}\,h^{(1)}_{3} +
    h^{(0)}_{1}\,h^{(2)}_{1}\,h^{(3)}_{3} + h^{(0)}_{1}\,h^{(1)}_{2}\,h^{(2)}_{2} +
    h^{(0)}_{1}\,h^{(1)}_{2}\,h^{(2)}_{3} + h^{(0)}_{1}\,h^{(1)}_{2}\,h^{(2)}_{4} \cr&\,\,\scriptstyle+
    h^{(0)}_{1}\,h^{(1)}_{2}\,h^{(3)}_{4} + h^{(0)}_{1}\,h^{(2)}_{2}\,h^{(3)}_{2} +
    h^{(0)}_{1}\,h^{(2)}_{2}\,h^{(1)}_{3} + h^{(0)}_{1}\,h^{(2)}_{2}\,h^{(3)}_{3} +
    h^{(0)}_{1}\,h^{(3)}_{2}\,h^{(2)}_{3} + h^{(0)}_{1}\,h^{(3)}_{2}\,h^{(2)}_{4} \cr&\,\,\scriptstyle+
    h^{(0)}_{1}\,h^{(1)}_{3}\,h^{(2)}_{3} + h^{(0)}_{1}\,h^{(1)}_{3}\,h^{(2)}_{4} +
    h^{(0)}_{1}\,h^{(1)}_{3}\,h^{(3)}_{4} + h^{(0)}_{1}\,h^{(2)}_{3}\,h^{(3)}_{3} +
    h^{(0)}_{1}\,h^{(2)}_{3}\,h^{(3)}_{4} + h^{(0)}_{1}\,h^{(3)}_{3}\,h^{(2)}_{4} \cr&\,\,\scriptstyle+
    h^{(0)}_{1}\,h^{(2)}_{4}\,h^{(3)}_{4} + h^{(1)}_{1}\,h^{(2)}_{1}\,h^{(0)}_{2} +
    h^{(1)}_{1}\,h^{(2)}_{1}\,h^{(3)}_{2} + h^{(1)}_{1}\,h^{(2)}_{1}\,h^{(0)}_{3} +
    h^{(1)}_{1}\,h^{(2)}_{1}\,h^{(0)}_{4} + h^{(1)}_{1}\,h^{(0)}_{2}\,h^{(2)}_{2} \cr&\,\,\scriptstyle+
    h^{(1)}_{1}\,h^{(0)}_{2}\,h^{(3)}_{2} + h^{(1)}_{1}\,h^{(0)}_{2}\,h^{(3)}_{3} +
    h^{(1)}_{1}\,h^{(0)}_{2}\,h^{(3)}_{4} + h^{(1)}_{1}\,h^{(2)}_{2}\,h^{(3)}_{2} +
    h^{(1)}_{1}\,h^{(2)}_{2}\,h^{(0)}_{3} + h^{(1)}_{1}\,h^{(2)}_{2}\,h^{(0)}_{4} \cr&\,\,\scriptstyle+
    h^{(1)}_{1}\,h^{(3)}_{2}\,h^{(2)}_{3} + h^{(1)}_{1}\,h^{(0)}_{3}\,h^{(2)}_{3} +
    h^{(1)}_{1}\,h^{(0)}_{3}\,h^{(3)}_{3} + h^{(1)}_{1}\,h^{(0)}_{3}\,h^{(3)}_{4} +
    h^{(1)}_{1}\,h^{(2)}_{3}\,h^{(3)}_{3} + h^{(1)}_{1}\,h^{(2)}_{3}\,h^{(0)}_{4} \cr&\,\,\scriptstyle+
    h^{(1)}_{1}\,h^{(2)}_{3}\,h^{(3)}_{4} + h^{(1)}_{1}\,h^{(3)}_{3}\,h^{(0)}_{4} +
    h^{(1)}_{1}\,h^{(0)}_{4}\,h^{(3)}_{4} + h^{(2)}_{1}\,h^{(0)}_{2}\,h^{(1)}_{2} +
    h^{(2)}_{1}\,h^{(0)}_{2}\,h^{(1)}_{3} + h^{(2)}_{1}\,h^{(0)}_{2}\,h^{(3)}_{3} \cr&\,\,\scriptstyle+
    h^{(2)}_{1}\,h^{(0)}_{2}\,h^{(1)}_{4} + h^{(2)}_{1}\,h^{(1)}_{2}\,h^{(3)}_{2} +
    h^{(2)}_{1}\,h^{(1)}_{2}\,h^{(0)}_{3} + h^{(2)}_{1}\,h^{(1)}_{2}\,h^{(0)}_{4} +
    h^{(2)}_{1}\,h^{(3)}_{2}\,h^{(1)}_{3} + h^{(2)}_{1}\,h^{(3)}_{2}\,h^{(1)}_{4} \cr&\,\,\scriptstyle+
    h^{(2)}_{1}\,h^{(0)}_{3}\,h^{(1)}_{3} + h^{(2)}_{1}\,h^{(0)}_{3}\,h^{(3)}_{3} +
    h^{(2)}_{1}\,h^{(0)}_{3}\,h^{(1)}_{4} + h^{(2)}_{1}\,h^{(1)}_{3}\,h^{(0)}_{4} +
    h^{(2)}_{1}\,h^{(3)}_{3}\,h^{(0)}_{4} + h^{(2)}_{1}\,h^{(3)}_{3}\,h^{(1)}_{4} \cr&\,\,\scriptstyle+
    h^{(2)}_{1}\,h^{(0)}_{4}\,h^{(1)}_{4} + h^{(0)}_{2}\,h^{(1)}_{2}\,h^{(2)}_{2} +
    h^{(0)}_{2}\,h^{(1)}_{2}\,h^{(3)}_{2} + h^{(0)}_{2}\,h^{(1)}_{2}\,h^{(2)}_{4} +
    h^{(0)}_{2}\,h^{(1)}_{2}\,h^{(3)}_{4} + h^{(0)}_{2}\,h^{(2)}_{2}\,h^{(1)}_{3} \cr&\,\,\scriptstyle+
    h^{(0)}_{2}\,h^{(2)}_{2}\,h^{(3)}_{3} + h^{(0)}_{2}\,h^{(2)}_{2}\,h^{(1)}_{4} +
    h^{(0)}_{2}\,h^{(3)}_{2}\,h^{(1)}_{3} + h^{(0)}_{2}\,h^{(3)}_{2}\,h^{(1)}_{4} +
    h^{(0)}_{2}\,h^{(3)}_{2}\,h^{(2)}_{4} + h^{(0)}_{2}\,h^{(1)}_{3}\,h^{(2)}_{4} \cr&\,\,\scriptstyle+
    h^{(0)}_{2}\,h^{(1)}_{3}\,h^{(3)}_{4} + h^{(0)}_{2}\,h^{(3)}_{3}\,h^{(1)}_{4} +
    h^{(0)}_{2}\,h^{(3)}_{3}\,h^{(2)}_{4} + h^{(0)}_{2}\,h^{(1)}_{4}\,h^{(3)}_{4} +
    h^{(0)}_{2}\,h^{(2)}_{4}\,h^{(3)}_{4} + h^{(1)}_{2}\,h^{(2)}_{2}\,h^{(3)}_{2} \cr&\,\,\scriptstyle+
    h^{(1)}_{2}\,h^{(2)}_{2}\,h^{(0)}_{3} + h^{(1)}_{2}\,h^{(2)}_{2}\,h^{(0)}_{4} +
    h^{(1)}_{2}\,h^{(3)}_{2}\,h^{(2)}_{3} + h^{(1)}_{2}\,h^{(3)}_{2}\,h^{(2)}_{4} +
    h^{(1)}_{2}\,h^{(0)}_{3}\,h^{(2)}_{3} + h^{(1)}_{2}\,h^{(0)}_{3}\,h^{(2)}_{4} \cr&\,\,\scriptstyle+
    h^{(1)}_{2}\,h^{(0)}_{3}\,h^{(3)}_{4} + h^{(1)}_{2}\,h^{(2)}_{3}\,h^{(0)}_{4} +
    h^{(1)}_{2}\,h^{(2)}_{3}\,h^{(3)}_{4} + h^{(1)}_{2}\,h^{(0)}_{4}\,h^{(2)}_{4} +
    h^{(1)}_{2}\,h^{(0)}_{4}\,h^{(3)}_{4} + h^{(2)}_{2}\,h^{(3)}_{2}\,h^{(1)}_{3} \cr&\,\,\scriptstyle+
    h^{(2)}_{2}\,h^{(3)}_{2}\,h^{(1)}_{4} + h^{(2)}_{2}\,h^{(0)}_{3}\,h^{(1)}_{3} +
    h^{(2)}_{2}\,h^{(0)}_{3}\,h^{(3)}_{3} + h^{(2)}_{2}\,h^{(0)}_{3}\,h^{(1)}_{4} +
    h^{(2)}_{2}\,h^{(1)}_{3}\,h^{(0)}_{4} + h^{(2)}_{2}\,h^{(3)}_{3}\,h^{(0)}_{4} \cr&\,\,\scriptstyle+
    h^{(2)}_{2}\,h^{(3)}_{3}\,h^{(1)}_{4} + h^{(2)}_{2}\,h^{(0)}_{4}\,h^{(1)}_{4} +
    h^{(3)}_{2}\,h^{(1)}_{3}\,h^{(2)}_{3} + h^{(3)}_{2}\,h^{(1)}_{3}\,h^{(2)}_{4} +
    h^{(3)}_{2}\,h^{(2)}_{3}\,h^{(1)}_{4} + h^{(3)}_{2}\,h^{(1)}_{4}\,h^{(2)}_{4} \cr&\,\,\scriptstyle+
    h^{(0)}_{3}\,h^{(1)}_{3}\,h^{(2)}_{3} + h^{(0)}_{3}\,h^{(1)}_{3}\,h^{(2)}_{4} +
    h^{(0)}_{3}\,h^{(1)}_{3}\,h^{(3)}_{4} + h^{(0)}_{3}\,h^{(2)}_{3}\,h^{(3)}_{3} +
    h^{(0)}_{3}\,h^{(2)}_{3}\,h^{(1)}_{4} + h^{(0)}_{3}\,h^{(3)}_{3}\,h^{(2)}_{4} \cr&\,\,\scriptstyle+
    h^{(0)}_{3}\,h^{(1)}_{4}\,h^{(2)}_{4} + h^{(0)}_{3}\,h^{(1)}_{4}\,h^{(3)}_{4} +
    h^{(0)}_{3}\,h^{(2)}_{4}\,h^{(3)}_{4} + h^{(1)}_{3}\,h^{(2)}_{3}\,h^{(0)}_{4} +
    h^{(1)}_{3}\,h^{(2)}_{3}\,h^{(3)}_{4} + h^{(1)}_{3}\,h^{(0)}_{4}\,h^{(2)}_{4} \cr&\,\,\scriptstyle+
    h^{(1)}_{3}\,h^{(0)}_{4}\,h^{(3)}_{4} + h^{(2)}_{3}\,h^{(3)}_{3}\,h^{(0)}_{4} +
    h^{(2)}_{3}\,h^{(3)}_{3}\,h^{(1)}_{4} + h^{(2)}_{3}\,h^{(0)}_{4}\,h^{(1)}_{4} +
    h^{(2)}_{3}\,h^{(1)}_{4}\,h^{(3)}_{4} + h^{(3)}_{3}\,h^{(0)}_{4}\,h^{(2)}_{4} \cr&\,\,\scriptstyle+
    h^{(3)}_{3}\,h^{(1)}_{4}\,h^{(2)}_{4} + h^{(0)}_{4}\,h^{(1)}_{4}\,h^{(2)}_{4} +
    h^{(0)}_{4}\,h^{(1)}_{4}\,h^{(3)}_{4} + h^{(0)}_{4}\,h^{(2)}_{4}\,h^{(3)}_{4} +
    h^{(1)}_{4}\,h^{(2)}_{4}\,h^{(3)}_{4}
\cr}$$
There are 125 terms in the above expression, all of which are positive.
Therefore $D_4 > 0$ since $h^{(j)}_{i}>0$ for all $i$ and $j$.

\vfill\break
$$\scriptstyle  A_3^\tau= \pmatrix
\scriptstyle h^{(0)}_{1} + h^{(0)}_{2} + h^{(0)}_{3} + h^{(0)}_{4} + h^{(1)}_{4} &
\scriptstyle -h^{(0)}_{1} - h^{(1)}_{1} - h^{(1)}_{2} - h^{(1)}_{3} - h^{(1)}_{4} & 
\scriptstyle h^{(1)}_{2} - h^{(0)}_{3} \cr
\scriptstyle h^{(0)}_{1} + h^{(0)}_{2} + h^{(0)}_{3} + h^{(2)}_{3} + h^{(0)}_{4} &
\scriptstyle -h^{(0)}_{1} + h^{(2)}_{4} &  h^{(2)}_{1} - h^{(0)}_{3} \cr
\scriptstyle h^{(0)}_{1} + h^{(0)}_{2} + h^{(3)}_{2} + h^{(0)}_{3} + h^{(0)}_{4} &
\scriptstyle -h^{(0)}_{1} + h^{(3)}_{3} &  
\scriptstyle -h^{(3)}_{1} - h^{(3)}_{2} - h^{(0)}_{3} - h^{(3)}_{3} - h^{(3)}_{4} \cr
\endpmatrix
$$

$$\eqalign{
\scriptstyle D_3&\scriptstyle = \det(A_3) \cr&= \scriptstyle
   -h^{(0)}_{1}\,h^{(1)}_{1}\,h^{(2)}_{1} - h^{(0)}_{1}\,h^{(1)}_{1}\,h^{(3)}_{1} -
    h^{(0)}_{1}\,h^{(1)}_{1}\,h^{(3)}_{2} - h^{(0)}_{1}\,h^{(1)}_{1}\,h^{(3)}_{3} -
    h^{(0)}_{1}\,h^{(1)}_{1}\,h^{(3)}_{4} - h^{(0)}_{1}\,h^{(2)}_{1}\,h^{(1)}_{2} \cr&\,\,\scriptstyle-
    h^{(0)}_{1}\,h^{(2)}_{1}\,h^{(3)}_{2} - h^{(0)}_{1}\,h^{(2)}_{1}\,h^{(1)}_{3} -
    h^{(0)}_{1}\,h^{(2)}_{1}\,h^{(3)}_{3} - h^{(0)}_{1}\,h^{(3)}_{1}\,h^{(1)}_{2} -
    h^{(0)}_{1}\,h^{(3)}_{1}\,h^{(1)}_{3} - h^{(0)}_{1}\,h^{(3)}_{1}\,h^{(2)}_{3} \cr&\,\,\scriptstyle-
    h^{(0)}_{1}\,h^{(3)}_{1}\,h^{(2)}_{4} - h^{(0)}_{1}\,h^{(1)}_{2}\,h^{(2)}_{3} -
    h^{(0)}_{1}\,h^{(1)}_{2}\,h^{(2)}_{4} - h^{(0)}_{1}\,h^{(1)}_{2}\,h^{(3)}_{4} -
    h^{(0)}_{1}\,h^{(3)}_{2}\,h^{(1)}_{3} - h^{(0)}_{1}\,h^{(3)}_{2}\,h^{(2)}_{3} \cr&\,\,\scriptstyle-
    h^{(0)}_{1}\,h^{(3)}_{2}\,h^{(2)}_{4} - h^{(0)}_{1}\,h^{(1)}_{3}\,h^{(3)}_{3} -
    h^{(0)}_{1}\,h^{(1)}_{3}\,h^{(3)}_{4} - h^{(0)}_{1}\,h^{(2)}_{3}\,h^{(3)}_{3} -
    h^{(0)}_{1}\,h^{(2)}_{3}\,h^{(3)}_{4} - h^{(0)}_{1}\,h^{(3)}_{3}\,h^{(2)}_{4} \cr&\,\,\scriptstyle-
    h^{(0)}_{1}\,h^{(2)}_{4}\,h^{(3)}_{4} - h^{(1)}_{1}\,h^{(2)}_{1}\,h^{(0)}_{2} -
    h^{(1)}_{1}\,h^{(2)}_{1}\,h^{(3)}_{2} - h^{(1)}_{1}\,h^{(2)}_{1}\,h^{(0)}_{3} -
    h^{(1)}_{1}\,h^{(2)}_{1}\,h^{(0)}_{4} - h^{(1)}_{1}\,h^{(3)}_{1}\,h^{(0)}_{2} \cr&\,\,\scriptstyle-
    h^{(1)}_{1}\,h^{(3)}_{1}\,h^{(0)}_{3} - h^{(1)}_{1}\,h^{(3)}_{1}\,h^{(2)}_{3} -
    h^{(1)}_{1}\,h^{(3)}_{1}\,h^{(0)}_{4} - h^{(1)}_{1}\,h^{(0)}_{2}\,h^{(3)}_{2} -
    h^{(1)}_{1}\,h^{(0)}_{2}\,h^{(3)}_{3} - h^{(1)}_{1}\,h^{(0)}_{2}\,h^{(3)}_{4} \cr&\,\,\scriptstyle-
    h^{(1)}_{1}\,h^{(3)}_{2}\,h^{(2)}_{3} - h^{(1)}_{1}\,h^{(3)}_{2}\,h^{(0)}_{4} -
    h^{(1)}_{1}\,h^{(0)}_{3}\,h^{(2)}_{3} - h^{(1)}_{1}\,h^{(0)}_{3}\,h^{(3)}_{3} -
    h^{(1)}_{1}\,h^{(0)}_{3}\,h^{(3)}_{4} - h^{(1)}_{1}\,h^{(2)}_{3}\,h^{(3)}_{3} \cr&\,\,\scriptstyle-
    h^{(1)}_{1}\,h^{(2)}_{3}\,h^{(3)}_{4} - h^{(1)}_{1}\,h^{(3)}_{3}\,h^{(0)}_{4} -
    h^{(1)}_{1}\,h^{(0)}_{4}\,h^{(3)}_{4} - h^{(2)}_{1}\,h^{(0)}_{2}\,h^{(1)}_{2} -
    h^{(2)}_{1}\,h^{(0)}_{2}\,h^{(1)}_{3} - h^{(2)}_{1}\,h^{(0)}_{2}\,h^{(3)}_{3} \cr&\,\,\scriptstyle-
    h^{(2)}_{1}\,h^{(0)}_{2}\,h^{(1)}_{4} - h^{(2)}_{1}\,h^{(1)}_{2}\,h^{(3)}_{2} -
    h^{(2)}_{1}\,h^{(1)}_{2}\,h^{(0)}_{3} - h^{(2)}_{1}\,h^{(1)}_{2}\,h^{(0)}_{4} -
    h^{(2)}_{1}\,h^{(3)}_{2}\,h^{(1)}_{3} - h^{(2)}_{1}\,h^{(3)}_{2}\,h^{(1)}_{4} \cr&\,\,\scriptstyle-
    h^{(2)}_{1}\,h^{(0)}_{3}\,h^{(1)}_{3} - h^{(2)}_{1}\,h^{(0)}_{3}\,h^{(3)}_{3} -
    h^{(2)}_{1}\,h^{(0)}_{3}\,h^{(1)}_{4} - h^{(2)}_{1}\,h^{(1)}_{3}\,h^{(0)}_{4} -
    h^{(2)}_{1}\,h^{(3)}_{3}\,h^{(0)}_{4} - h^{(2)}_{1}\,h^{(3)}_{3}\,h^{(1)}_{4} \cr&\,\,\scriptstyle-
    h^{(2)}_{1}\,h^{(0)}_{4}\,h^{(1)}_{4} - h^{(3)}_{1}\,h^{(0)}_{2}\,h^{(1)}_{2} -
    h^{(3)}_{1}\,h^{(0)}_{2}\,h^{(1)}_{3} - h^{(3)}_{1}\,h^{(0)}_{2}\,h^{(1)}_{4} -
    h^{(3)}_{1}\,h^{(0)}_{2}\,h^{(2)}_{4} - h^{(3)}_{1}\,h^{(1)}_{2}\,h^{(0)}_{3} \cr&\,\,\scriptstyle-
    h^{(3)}_{1}\,h^{(1)}_{2}\,h^{(2)}_{3} - h^{(3)}_{1}\,h^{(1)}_{2}\,h^{(0)}_{4} -
    h^{(3)}_{1}\,h^{(0)}_{3}\,h^{(1)}_{3} - h^{(3)}_{1}\,h^{(0)}_{3}\,h^{(1)}_{4} -
    h^{(3)}_{1}\,h^{(0)}_{3}\,h^{(2)}_{4} - h^{(3)}_{1}\,h^{(1)}_{3}\,h^{(2)}_{3} \cr&\,\,\scriptstyle-
    h^{(3)}_{1}\,h^{(1)}_{3}\,h^{(0)}_{4} - h^{(3)}_{1}\,h^{(2)}_{3}\,h^{(1)}_{4} -
    h^{(3)}_{1}\,h^{(0)}_{4}\,h^{(1)}_{4} - h^{(3)}_{1}\,h^{(0)}_{4}\,h^{(2)}_{4} -
    h^{(3)}_{1}\,h^{(1)}_{4}\,h^{(2)}_{4} - h^{(0)}_{2}\,h^{(1)}_{2}\,h^{(3)}_{2} \cr&\,\,\scriptstyle-
    h^{(0)}_{2}\,h^{(1)}_{2}\,h^{(2)}_{4} - h^{(0)}_{2}\,h^{(1)}_{2}\,h^{(3)}_{4} -
    h^{(0)}_{2}\,h^{(3)}_{2}\,h^{(1)}_{3} - h^{(0)}_{2}\,h^{(3)}_{2}\,h^{(1)}_{4} -
    h^{(0)}_{2}\,h^{(3)}_{2}\,h^{(2)}_{4} - h^{(0)}_{2}\,h^{(1)}_{3}\,h^{(3)}_{3} \cr&\,\,\scriptstyle-
    h^{(0)}_{2}\,h^{(1)}_{3}\,h^{(3)}_{4} - h^{(0)}_{2}\,h^{(3)}_{3}\,h^{(1)}_{4} -
    h^{(0)}_{2}\,h^{(3)}_{3}\,h^{(2)}_{4} - h^{(0)}_{2}\,h^{(1)}_{4}\,h^{(3)}_{4} -
    h^{(0)}_{2}\,h^{(2)}_{4}\,h^{(3)}_{4} - h^{(1)}_{2}\,h^{(3)}_{2}\,h^{(2)}_{3} \cr&\,\,\scriptstyle-
    h^{(1)}_{2}\,h^{(3)}_{2}\,h^{(0)}_{4} - h^{(1)}_{2}\,h^{(3)}_{2}\,h^{(2)}_{4} -
    h^{(1)}_{2}\,h^{(0)}_{3}\,h^{(2)}_{3} - h^{(1)}_{2}\,h^{(0)}_{3}\,h^{(2)}_{4} -
    h^{(1)}_{2}\,h^{(0)}_{3}\,h^{(3)}_{4} - h^{(1)}_{2}\,h^{(2)}_{3}\,h^{(3)}_{4} \cr&\,\,\scriptstyle-
    h^{(1)}_{2}\,h^{(0)}_{4}\,h^{(2)}_{4} - h^{(1)}_{2}\,h^{(0)}_{4}\,h^{(3)}_{4} -
    h^{(3)}_{2}\,h^{(1)}_{3}\,h^{(2)}_{3} - h^{(3)}_{2}\,h^{(1)}_{3}\,h^{(0)}_{4} -
    h^{(3)}_{2}\,h^{(2)}_{3}\,h^{(1)}_{4} - h^{(3)}_{2}\,h^{(0)}_{4}\,h^{(1)}_{4} \cr&\,\,\scriptstyle-
    h^{(3)}_{2}\,h^{(0)}_{4}\,h^{(2)}_{4} - h^{(3)}_{2}\,h^{(1)}_{4}\,h^{(2)}_{4} -
    h^{(0)}_{3}\,h^{(1)}_{3}\,h^{(2)}_{3} - h^{(0)}_{3}\,h^{(1)}_{3}\,h^{(3)}_{3} -
    h^{(0)}_{3}\,h^{(1)}_{3}\,h^{(3)}_{4} - h^{(0)}_{3}\,h^{(2)}_{3}\,h^{(3)}_{3} \cr&\,\,\scriptstyle-
    h^{(0)}_{3}\,h^{(2)}_{3}\,h^{(1)}_{4} - h^{(0)}_{3}\,h^{(3)}_{3}\,h^{(2)}_{4} -
    h^{(0)}_{3}\,h^{(1)}_{4}\,h^{(2)}_{4} - h^{(0)}_{3}\,h^{(1)}_{4}\,h^{(3)}_{4} -
    h^{(0)}_{3}\,h^{(2)}_{4}\,h^{(3)}_{4} - h^{(1)}_{3}\,h^{(2)}_{3}\,h^{(3)}_{3} \cr&\,\,\scriptstyle-
    h^{(1)}_{3}\,h^{(2)}_{3}\,h^{(3)}_{4} - h^{(1)}_{3}\,h^{(3)}_{3}\,h^{(0)}_{4} -
    h^{(1)}_{3}\,h^{(0)}_{4}\,h^{(3)}_{4} - h^{(2)}_{3}\,h^{(3)}_{3}\,h^{(1)}_{4} -
    h^{(2)}_{3}\,h^{(1)}_{4}\,h^{(3)}_{4} - h^{(3)}_{3}\,h^{(0)}_{4}\,h^{(1)}_{4} \cr&\,\,\scriptstyle-
    h^{(3)}_{3}\,h^{(0)}_{4}\,h^{(2)}_{4} - h^{(3)}_{3}\,h^{(1)}_{4}\,h^{(2)}_{4} -
    h^{(0)}_{4}\,h^{(1)}_{4}\,h^{(3)}_{4} - h^{(0)}_{4}\,h^{(2)}_{4}\,h^{(3)}_{4} -
    h^{(1)}_{4}\,h^{(2)}_{4}\,h^{(3)}_{4}
\cr}$$
There are 125 terms in the above expression, all of which are negative.
Therefore $D_3 < 0$ since $h^{(j)}_{i}>0$ for all $i$ and $j$.

\vfill\break
$$\scriptstyle  A_2^\tau= \pmatrix
\scriptstyle h^{(0)}_{1} + h^{(0)}_{2} + h^{(0)}_{3} + h^{(0)}_{4} + h^{(1)}_{4} &
\scriptstyle h^{(1)}_{1} - h^{(0)}_{2}  & 
\scriptstyle h^{(1)}_{2} - h^{(0)}_{3} \cr
\scriptstyle h^{(0)}_{1} + h^{(0)}_{2} + h^{(0)}_{3} + h^{(2)}_{3} + h^{(0)}_{4} &
\scriptstyle -h^{(2)}_{1} - h^{(0)}_{2} - h^{(2)}_{2} - h^{(2)}_{3} - h^{(2)}_{4} &
\scriptstyle h^{(2)}_{1} - h^{(0)}_{3} \cr
\scriptstyle h^{(0)}_{1} + h^{(0)}_{2} + h^{(3)}_{2} + h^{(0)}_{3} + h^{(0)}_{4} &
\scriptstyle -h^{(0)}_{2} + h^{(3)}_{4} &  
\scriptstyle -h^{(3)}_{1} - h^{(3)}_{2} - h^{(0)}_{3} - h^{(3)}_{3} - h^{(3)}_{4} \cr
\endpmatrix
$$

$$\eqalign{
\scriptstyle D_2&\scriptstyle = \det(A_2) \cr&= \scriptstyle
    h^{(0)}_{1}\,h^{(1)}_{1}\,h^{(2)}_{1} + h^{(0)}_{1}\,h^{(1)}_{1}\,h^{(3)}_{1} +
    h^{(0)}_{1}\,h^{(1)}_{1}\,h^{(3)}_{2} + h^{(0)}_{1}\,h^{(1)}_{1}\,h^{(3)}_{3} +
    h^{(0)}_{1}\,h^{(1)}_{1}\,h^{(3)}_{4} + h^{(0)}_{1}\,h^{(2)}_{1}\,h^{(3)}_{1} \cr&\,\,\scriptstyle+
    h^{(0)}_{1}\,h^{(2)}_{1}\,h^{(1)}_{2} + h^{(0)}_{1}\,h^{(2)}_{1}\,h^{(3)}_{2} +
    h^{(0)}_{1}\,h^{(2)}_{1}\,h^{(3)}_{3} + h^{(0)}_{1}\,h^{(3)}_{1}\,h^{(2)}_{2} +
    h^{(0)}_{1}\,h^{(3)}_{1}\,h^{(2)}_{3} + h^{(0)}_{1}\,h^{(3)}_{1}\,h^{(2)}_{4} \cr&\,\,\scriptstyle+
    h^{(0)}_{1}\,h^{(1)}_{2}\,h^{(2)}_{2} + h^{(0)}_{1}\,h^{(1)}_{2}\,h^{(2)}_{3} +
    h^{(0)}_{1}\,h^{(1)}_{2}\,h^{(2)}_{4} + h^{(0)}_{1}\,h^{(1)}_{2}\,h^{(3)}_{4} +
    h^{(0)}_{1}\,h^{(2)}_{2}\,h^{(3)}_{2} + h^{(0)}_{1}\,h^{(2)}_{2}\,h^{(3)}_{3} \cr&\,\,\scriptstyle+
    h^{(0)}_{1}\,h^{(2)}_{2}\,h^{(3)}_{4} + h^{(0)}_{1}\,h^{(3)}_{2}\,h^{(2)}_{3} +
    h^{(0)}_{1}\,h^{(3)}_{2}\,h^{(2)}_{4} + h^{(0)}_{1}\,h^{(2)}_{3}\,h^{(3)}_{3} +
    h^{(0)}_{1}\,h^{(2)}_{3}\,h^{(3)}_{4} + h^{(0)}_{1}\,h^{(3)}_{3}\,h^{(2)}_{4} \cr&\,\,\scriptstyle+
    h^{(0)}_{1}\,h^{(2)}_{4}\,h^{(3)}_{4} + h^{(1)}_{1}\,h^{(2)}_{1}\,h^{(0)}_{2} +
    h^{(1)}_{1}\,h^{(2)}_{1}\,h^{(3)}_{2} + h^{(1)}_{1}\,h^{(2)}_{1}\,h^{(0)}_{3} +
    h^{(1)}_{1}\,h^{(2)}_{1}\,h^{(0)}_{4} + h^{(1)}_{1}\,h^{(3)}_{1}\,h^{(0)}_{2} \cr&\,\,\scriptstyle+
    h^{(1)}_{1}\,h^{(3)}_{1}\,h^{(0)}_{3} + h^{(1)}_{1}\,h^{(3)}_{1}\,h^{(2)}_{3} +
    h^{(1)}_{1}\,h^{(3)}_{1}\,h^{(0)}_{4} + h^{(1)}_{1}\,h^{(0)}_{2}\,h^{(3)}_{2} +
    h^{(1)}_{1}\,h^{(0)}_{2}\,h^{(3)}_{3} + h^{(1)}_{1}\,h^{(0)}_{2}\,h^{(3)}_{4} \cr&\,\,\scriptstyle+
    h^{(1)}_{1}\,h^{(3)}_{2}\,h^{(2)}_{3} + h^{(1)}_{1}\,h^{(3)}_{2}\,h^{(0)}_{4} +
    h^{(1)}_{1}\,h^{(0)}_{3}\,h^{(2)}_{3} + h^{(1)}_{1}\,h^{(0)}_{3}\,h^{(3)}_{3} +
    h^{(1)}_{1}\,h^{(0)}_{3}\,h^{(3)}_{4} + h^{(1)}_{1}\,h^{(2)}_{3}\,h^{(3)}_{3} \cr&\,\,\scriptstyle+
    h^{(1)}_{1}\,h^{(2)}_{3}\,h^{(3)}_{4} + h^{(1)}_{1}\,h^{(3)}_{3}\,h^{(0)}_{4} +
    h^{(1)}_{1}\,h^{(0)}_{4}\,h^{(3)}_{4} + h^{(2)}_{1}\,h^{(3)}_{1}\,h^{(0)}_{2} +
    h^{(2)}_{1}\,h^{(3)}_{1}\,h^{(0)}_{3} + h^{(2)}_{1}\,h^{(3)}_{1}\,h^{(0)}_{4} \cr&\,\,\scriptstyle+
    h^{(2)}_{1}\,h^{(3)}_{1}\,h^{(1)}_{4} + h^{(2)}_{1}\,h^{(0)}_{2}\,h^{(1)}_{2} +
    h^{(2)}_{1}\,h^{(0)}_{2}\,h^{(3)}_{3} + h^{(2)}_{1}\,h^{(0)}_{2}\,h^{(1)}_{4} +
    h^{(2)}_{1}\,h^{(1)}_{2}\,h^{(3)}_{2} + h^{(2)}_{1}\,h^{(1)}_{2}\,h^{(0)}_{3} \cr&\,\,\scriptstyle+
    h^{(2)}_{1}\,h^{(1)}_{2}\,h^{(0)}_{4} + h^{(2)}_{1}\,h^{(3)}_{2}\,h^{(0)}_{4} +
    h^{(2)}_{1}\,h^{(3)}_{2}\,h^{(1)}_{4} + h^{(2)}_{1}\,h^{(0)}_{3}\,h^{(3)}_{3} +
    h^{(2)}_{1}\,h^{(0)}_{3}\,h^{(1)}_{4} + h^{(2)}_{1}\,h^{(3)}_{3}\,h^{(0)}_{4} \cr&\,\,\scriptstyle+
    h^{(2)}_{1}\,h^{(3)}_{3}\,h^{(1)}_{4} + h^{(3)}_{1}\,h^{(0)}_{2}\,h^{(2)}_{2} +
    h^{(3)}_{1}\,h^{(0)}_{2}\,h^{(1)}_{4} + h^{(3)}_{1}\,h^{(0)}_{2}\,h^{(2)}_{4} +
    h^{(3)}_{1}\,h^{(2)}_{2}\,h^{(0)}_{3} + h^{(3)}_{1}\,h^{(2)}_{2}\,h^{(0)}_{4} \cr&\,\,\scriptstyle+
    h^{(3)}_{1}\,h^{(2)}_{2}\,h^{(1)}_{4} + h^{(3)}_{1}\,h^{(0)}_{3}\,h^{(2)}_{3} +
    h^{(3)}_{1}\,h^{(0)}_{3}\,h^{(2)}_{4} + h^{(3)}_{1}\,h^{(2)}_{3}\,h^{(0)}_{4} +
    h^{(3)}_{1}\,h^{(2)}_{3}\,h^{(1)}_{4} + h^{(3)}_{1}\,h^{(0)}_{4}\,h^{(2)}_{4} \cr&\,\,\scriptstyle+
    h^{(3)}_{1}\,h^{(1)}_{4}\,h^{(2)}_{4} + h^{(0)}_{2}\,h^{(1)}_{2}\,h^{(2)}_{2} +
    h^{(0)}_{2}\,h^{(1)}_{2}\,h^{(3)}_{2} + h^{(0)}_{2}\,h^{(1)}_{2}\,h^{(2)}_{4} +
    h^{(0)}_{2}\,h^{(1)}_{2}\,h^{(3)}_{4} + h^{(0)}_{2}\,h^{(2)}_{2}\,h^{(3)}_{2} \cr&\,\,\scriptstyle+
    h^{(0)}_{2}\,h^{(2)}_{2}\,h^{(3)}_{3} + h^{(0)}_{2}\,h^{(2)}_{2}\,h^{(3)}_{4} +
    h^{(0)}_{2}\,h^{(3)}_{2}\,h^{(1)}_{4} + h^{(0)}_{2}\,h^{(3)}_{2}\,h^{(2)}_{4} +
    h^{(0)}_{2}\,h^{(3)}_{3}\,h^{(1)}_{4} + h^{(0)}_{2}\,h^{(3)}_{3}\,h^{(2)}_{4} \cr&\,\,\scriptstyle+
    h^{(0)}_{2}\,h^{(1)}_{4}\,h^{(3)}_{4} + h^{(0)}_{2}\,h^{(2)}_{4}\,h^{(3)}_{4} +
    h^{(1)}_{2}\,h^{(2)}_{2}\,h^{(3)}_{2} + h^{(1)}_{2}\,h^{(2)}_{2}\,h^{(0)}_{3} +
    h^{(1)}_{2}\,h^{(2)}_{2}\,h^{(0)}_{4} + h^{(1)}_{2}\,h^{(3)}_{2}\,h^{(2)}_{3} \cr&\,\,\scriptstyle+
    h^{(1)}_{2}\,h^{(3)}_{2}\,h^{(2)}_{4} + h^{(1)}_{2}\,h^{(0)}_{3}\,h^{(2)}_{3} +
    h^{(1)}_{2}\,h^{(0)}_{3}\,h^{(2)}_{4} + h^{(1)}_{2}\,h^{(0)}_{3}\,h^{(3)}_{4} +
    h^{(1)}_{2}\,h^{(2)}_{3}\,h^{(0)}_{4} + h^{(1)}_{2}\,h^{(2)}_{3}\,h^{(3)}_{4} \cr&\,\,\scriptstyle+
    h^{(1)}_{2}\,h^{(0)}_{4}\,h^{(2)}_{4} + h^{(1)}_{2}\,h^{(0)}_{4}\,h^{(3)}_{4} +
    h^{(2)}_{2}\,h^{(3)}_{2}\,h^{(0)}_{4} + h^{(2)}_{2}\,h^{(3)}_{2}\,h^{(1)}_{4} +
    h^{(2)}_{2}\,h^{(0)}_{3}\,h^{(3)}_{3} + h^{(2)}_{2}\,h^{(0)}_{3}\,h^{(1)}_{4} \cr&\,\,\scriptstyle+
    h^{(2)}_{2}\,h^{(0)}_{3}\,h^{(3)}_{4} + h^{(2)}_{2}\,h^{(3)}_{3}\,h^{(0)}_{4} +
    h^{(2)}_{2}\,h^{(3)}_{3}\,h^{(1)}_{4} + h^{(2)}_{2}\,h^{(0)}_{4}\,h^{(3)}_{4} +
    h^{(2)}_{2}\,h^{(1)}_{4}\,h^{(3)}_{4} + h^{(3)}_{2}\,h^{(2)}_{3}\,h^{(0)}_{4} \cr&\,\,\scriptstyle+
    h^{(3)}_{2}\,h^{(2)}_{3}\,h^{(1)}_{4} + h^{(3)}_{2}\,h^{(0)}_{4}\,h^{(2)}_{4} +
    h^{(3)}_{2}\,h^{(1)}_{4}\,h^{(2)}_{4} + h^{(0)}_{3}\,h^{(2)}_{3}\,h^{(3)}_{3} +
    h^{(0)}_{3}\,h^{(2)}_{3}\,h^{(1)}_{4} + h^{(0)}_{3}\,h^{(3)}_{3}\,h^{(2)}_{4} \cr&\,\,\scriptstyle+
    h^{(0)}_{3}\,h^{(1)}_{4}\,h^{(2)}_{4} + h^{(0)}_{3}\,h^{(1)}_{4}\,h^{(3)}_{4} +
    h^{(0)}_{3}\,h^{(2)}_{4}\,h^{(3)}_{4} + h^{(2)}_{3}\,h^{(3)}_{3}\,h^{(0)}_{4} +
    h^{(2)}_{3}\,h^{(3)}_{3}\,h^{(1)}_{4} + h^{(2)}_{3}\,h^{(0)}_{4}\,h^{(3)}_{4} \cr&\,\,\scriptstyle+
    h^{(2)}_{3}\,h^{(1)}_{4}\,h^{(3)}_{4} + h^{(3)}_{3}\,h^{(0)}_{4}\,h^{(2)}_{4} +
    h^{(3)}_{3}\,h^{(1)}_{4}\,h^{(2)}_{4} + h^{(0)}_{4}\,h^{(2)}_{4}\,h^{(3)}_{4} +
    h^{(1)}_{4}\,h^{(2)}_{4}\,h^{(3)}_{4}
\cr}$$
There are 125 terms in the above expression, all of which are positive.
Therefore $D_2 > 0$ since $h^{(j)}_{i}>0$ for all $i$ and $j$.

\vfill\break
$$\scriptstyle  A_1^\tau=\! \pmatrix
\scriptstyle -h^{(0)}_{1} - h^{(1)}_{1} - h^{(1)}_{2} - h^{(1)}_{3} - h^{(1)}_{4} &
\scriptstyle -h^{(1)}_{1} - h^{(0)}_{2}  & 
\scriptstyle h^{(1)}_{2} - h^{(0)}_{3} \cr
\scriptstyle -h^{(0)}_{1} + h^{(2)}_{4} &  
\scriptstyle -h^{(2)}_{1} - h^{(0)}_{2} - h^{(2)}_{2} -
\scriptstyle h^{(2)}_{3} - h^{(2)}_{4}  & 
\scriptstyle h^{(2)}_{1} - h^{(0)}_{3} \cr
\scriptstyle -h^{(0)}_{1} + h^{(3)}_{3} &  
\scriptstyle -h^{(0)}_{2} + h^{(3)}_{4} &  
\scriptstyle -h^{(3)}_{1} - h^{(3)}_{2} - h^{(0)}_{3} - h^{(3)}_{3} - h^{(3)}_{4} \cr
\endpmatrix
$$

$$\eqalign{
\scriptstyle D_1&\scriptstyle = \det(A_1) \cr&= \scriptstyle 
   -h^{(0)}_{1}\,h^{(1)}_{1}\,h^{(2)}_{1} - h^{(0)}_{1}\,h^{(1)}_{1}\,h^{(3)}_{1} -
    h^{(0)}_{1}\,h^{(1)}_{1}\,h^{(3)}_{2} - h^{(0)}_{1}\,h^{(1)}_{1}\,h^{(3)}_{3} -
    h^{(0)}_{1}\,h^{(1)}_{1}\,h^{(3)}_{4} - h^{(0)}_{1}\,h^{(2)}_{1}\,h^{(3)}_{1} \cr&\,\,\scriptstyle-
    h^{(0)}_{1}\,h^{(2)}_{1}\,h^{(1)}_{2} - h^{(0)}_{1}\,h^{(2)}_{1}\,h^{(3)}_{2} -
    h^{(0)}_{1}\,h^{(2)}_{1}\,h^{(3)}_{3} - h^{(0)}_{1}\,h^{(3)}_{1}\,h^{(2)}_{2} -
    h^{(0)}_{1}\,h^{(3)}_{1}\,h^{(2)}_{3} - h^{(0)}_{1}\,h^{(3)}_{1}\,h^{(2)}_{4} \cr&\,\,\scriptstyle-
    h^{(0)}_{1}\,h^{(1)}_{2}\,h^{(2)}_{2} - h^{(0)}_{1}\,h^{(1)}_{2}\,h^{(2)}_{3} -
    h^{(0)}_{1}\,h^{(1)}_{2}\,h^{(2)}_{4} - h^{(0)}_{1}\,h^{(1)}_{2}\,h^{(3)}_{4} -
    h^{(0)}_{1}\,h^{(2)}_{2}\,h^{(3)}_{2} - h^{(0)}_{1}\,h^{(2)}_{2}\,h^{(3)}_{3} \cr&\,\,\scriptstyle-
    h^{(0)}_{1}\,h^{(2)}_{2}\,h^{(3)}_{4} - h^{(0)}_{1}\,h^{(3)}_{2}\,h^{(2)}_{3} -
    h^{(0)}_{1}\,h^{(3)}_{2}\,h^{(2)}_{4} - h^{(0)}_{1}\,h^{(2)}_{3}\,h^{(3)}_{3} -
    h^{(0)}_{1}\,h^{(2)}_{3}\,h^{(3)}_{4} - h^{(0)}_{1}\,h^{(3)}_{3}\,h^{(2)}_{4} \cr&\,\,\scriptstyle-
    h^{(0)}_{1}\,h^{(2)}_{4}\,h^{(3)}_{4} - h^{(1)}_{1}\,h^{(2)}_{1}\,h^{(3)}_{1} -
    h^{(1)}_{1}\,h^{(2)}_{1}\,h^{(0)}_{2} - h^{(1)}_{1}\,h^{(2)}_{1}\,h^{(3)}_{2} -
    h^{(1)}_{1}\,h^{(2)}_{1}\,h^{(0)}_{3} - h^{(1)}_{1}\,h^{(3)}_{1}\,h^{(0)}_{2} \cr&\,\,\scriptstyle-
    h^{(1)}_{1}\,h^{(3)}_{1}\,h^{(2)}_{2} - h^{(1)}_{1}\,h^{(3)}_{1}\,h^{(2)}_{3} -
    h^{(1)}_{1}\,h^{(0)}_{2}\,h^{(3)}_{2} - h^{(1)}_{1}\,h^{(0)}_{2}\,h^{(3)}_{3} -
    h^{(1)}_{1}\,h^{(0)}_{2}\,h^{(3)}_{4} - h^{(1)}_{1}\,h^{(2)}_{2}\,h^{(3)}_{2} \cr&\,\,\scriptstyle-
    h^{(1)}_{1}\,h^{(2)}_{2}\,h^{(0)}_{3} - h^{(1)}_{1}\,h^{(2)}_{2}\,h^{(3)}_{3} -
    h^{(1)}_{1}\,h^{(2)}_{2}\,h^{(3)}_{4} - h^{(1)}_{1}\,h^{(3)}_{2}\,h^{(2)}_{3} -
    h^{(1)}_{1}\,h^{(0)}_{3}\,h^{(2)}_{3} - h^{(1)}_{1}\,h^{(0)}_{3}\,h^{(3)}_{3} \cr&\,\,\scriptstyle-
    h^{(1)}_{1}\,h^{(0)}_{3}\,h^{(3)}_{4} - h^{(1)}_{1}\,h^{(2)}_{3}\,h^{(3)}_{3} -
    h^{(1)}_{1}\,h^{(2)}_{3}\,h^{(3)}_{4} - h^{(2)}_{1}\,h^{(3)}_{1}\,h^{(1)}_{2} -
    h^{(2)}_{1}\,h^{(3)}_{1}\,h^{(1)}_{3} - h^{(2)}_{1}\,h^{(3)}_{1}\,h^{(1)}_{4} \cr&\,\,\scriptstyle-
    h^{(2)}_{1}\,h^{(0)}_{2}\,h^{(1)}_{2} - h^{(2)}_{1}\,h^{(0)}_{2}\,h^{(1)}_{3} -
    h^{(2)}_{1}\,h^{(0)}_{2}\,h^{(3)}_{3} - h^{(2)}_{1}\,h^{(0)}_{2}\,h^{(1)}_{4} -
    h^{(2)}_{1}\,h^{(1)}_{2}\,h^{(3)}_{2} - h^{(2)}_{1}\,h^{(1)}_{2}\,h^{(0)}_{3} \cr&\,\,\scriptstyle-
    h^{(2)}_{1}\,h^{(3)}_{2}\,h^{(1)}_{3} - h^{(2)}_{1}\,h^{(3)}_{2}\,h^{(1)}_{4} -
    h^{(2)}_{1}\,h^{(0)}_{3}\,h^{(1)}_{3} - h^{(2)}_{1}\,h^{(0)}_{3}\,h^{(3)}_{3} -
    h^{(2)}_{1}\,h^{(0)}_{3}\,h^{(1)}_{4} - h^{(2)}_{1}\,h^{(1)}_{3}\,h^{(3)}_{3} \cr&\,\,\scriptstyle-
    h^{(2)}_{1}\,h^{(3)}_{3}\,h^{(1)}_{4} - h^{(3)}_{1}\,h^{(0)}_{2}\,h^{(1)}_{2} -
    h^{(3)}_{1}\,h^{(0)}_{2}\,h^{(1)}_{3} - h^{(3)}_{1}\,h^{(0)}_{2}\,h^{(1)}_{4} -
    h^{(3)}_{1}\,h^{(0)}_{2}\,h^{(2)}_{4} - h^{(3)}_{1}\,h^{(1)}_{2}\,h^{(2)}_{2} \cr&\,\,\scriptstyle-
    h^{(3)}_{1}\,h^{(1)}_{2}\,h^{(2)}_{3} - h^{(3)}_{1}\,h^{(1)}_{2}\,h^{(2)}_{4} -
    h^{(3)}_{1}\,h^{(2)}_{2}\,h^{(1)}_{3} - h^{(3)}_{1}\,h^{(2)}_{2}\,h^{(1)}_{4} -
    h^{(3)}_{1}\,h^{(1)}_{3}\,h^{(2)}_{3} - h^{(3)}_{1}\,h^{(1)}_{3}\,h^{(2)}_{4} \cr&\,\,\scriptstyle-
    h^{(3)}_{1}\,h^{(2)}_{3}\,h^{(1)}_{4} - h^{(3)}_{1}\,h^{(1)}_{4}\,h^{(2)}_{4} -
    h^{(0)}_{2}\,h^{(1)}_{2}\,h^{(3)}_{2} - h^{(0)}_{2}\,h^{(1)}_{2}\,h^{(2)}_{4} -
    h^{(0)}_{2}\,h^{(1)}_{2}\,h^{(3)}_{4} - h^{(0)}_{2}\,h^{(3)}_{2}\,h^{(1)}_{3} \cr&\,\,\scriptstyle-
    h^{(0)}_{2}\,h^{(3)}_{2}\,h^{(1)}_{4} - h^{(0)}_{2}\,h^{(3)}_{2}\,h^{(2)}_{4} -
    h^{(0)}_{2}\,h^{(1)}_{3}\,h^{(3)}_{3} - h^{(0)}_{2}\,h^{(1)}_{3}\,h^{(3)}_{4} -
    h^{(0)}_{2}\,h^{(3)}_{3}\,h^{(1)}_{4} - h^{(0)}_{2}\,h^{(3)}_{3}\,h^{(2)}_{4} \cr&\,\,\scriptstyle-
    h^{(0)}_{2}\,h^{(1)}_{4}\,h^{(3)}_{4} - h^{(0)}_{2}\,h^{(2)}_{4}\,h^{(3)}_{4} -
    h^{(1)}_{2}\,h^{(2)}_{2}\,h^{(3)}_{2} - h^{(1)}_{2}\,h^{(2)}_{2}\,h^{(0)}_{3} -
    h^{(1)}_{2}\,h^{(2)}_{2}\,h^{(3)}_{4} - h^{(1)}_{2}\,h^{(3)}_{2}\,h^{(2)}_{3} \cr&\,\,\scriptstyle-
    h^{(1)}_{2}\,h^{(3)}_{2}\,h^{(2)}_{4} - h^{(1)}_{2}\,h^{(0)}_{3}\,h^{(2)}_{3} -
    h^{(1)}_{2}\,h^{(0)}_{3}\,h^{(2)}_{4} - h^{(1)}_{2}\,h^{(0)}_{3}\,h^{(3)}_{4} -
    h^{(1)}_{2}\,h^{(2)}_{3}\,h^{(3)}_{4} - h^{(2)}_{2}\,h^{(3)}_{2}\,h^{(1)}_{3} \cr&\,\,\scriptstyle-
    h^{(2)}_{2}\,h^{(3)}_{2}\,h^{(1)}_{4} - h^{(2)}_{2}\,h^{(0)}_{3}\,h^{(1)}_{3} -
    h^{(2)}_{2}\,h^{(0)}_{3}\,h^{(3)}_{3} - h^{(2)}_{2}\,h^{(0)}_{3}\,h^{(1)}_{4} -
    h^{(2)}_{2}\,h^{(1)}_{3}\,h^{(3)}_{3} - h^{(2)}_{2}\,h^{(1)}_{3}\,h^{(3)}_{4} \cr&\,\,\scriptstyle-
    h^{(2)}_{2}\,h^{(3)}_{3}\,h^{(1)}_{4} - h^{(2)}_{2}\,h^{(1)}_{4}\,h^{(3)}_{4} -
    h^{(3)}_{2}\,h^{(1)}_{3}\,h^{(2)}_{3} - h^{(3)}_{2}\,h^{(1)}_{3}\,h^{(2)}_{4} -
    h^{(3)}_{2}\,h^{(2)}_{3}\,h^{(1)}_{4} - h^{(3)}_{2}\,h^{(1)}_{4}\,h^{(2)}_{4} \cr&\,\,\scriptstyle-
    h^{(0)}_{3}\,h^{(1)}_{3}\,h^{(2)}_{3} - h^{(0)}_{3}\,h^{(1)}_{3}\,h^{(2)}_{4} -
    h^{(0)}_{3}\,h^{(1)}_{3}\,h^{(3)}_{4} - h^{(0)}_{3}\,h^{(2)}_{3}\,h^{(3)}_{3} -
    h^{(0)}_{3}\,h^{(2)}_{3}\,h^{(1)}_{4} - h^{(0)}_{3}\,h^{(3)}_{3}\,h^{(2)}_{4} \cr&\,\,\scriptstyle-
    h^{(0)}_{3}\,h^{(1)}_{4}\,h^{(2)}_{4} - h^{(0)}_{3}\,h^{(1)}_{4}\,h^{(3)}_{4} -
    h^{(0)}_{3}\,h^{(2)}_{4}\,h^{(3)}_{4} - h^{(1)}_{3}\,h^{(2)}_{3}\,h^{(3)}_{3} -
    h^{(1)}_{3}\,h^{(2)}_{3}\,h^{(3)}_{4} - h^{(1)}_{3}\,h^{(3)}_{3}\,h^{(2)}_{4} \cr&\,\,\scriptstyle-
    h^{(1)}_{3}\,h^{(2)}_{4}\,h^{(3)}_{4} - h^{(2)}_{3}\,h^{(3)}_{3}\,h^{(1)}_{4} -
    h^{(2)}_{3}\,h^{(1)}_{4}\,h^{(3)}_{4} - h^{(3)}_{3}\,h^{(1)}_{4}\,h^{(2)}_{4} -
    h^{(1)}_{4}\,h^{(2)}_{4}\,h^{(3)}_{4}
\cr}$$
There are 125 terms in the above expression, all of which are negative.
Therefore $D_1 < 0$ since $h^{(j)}_{i}>0$ for all $i$ and $j$.

\vfill\break
Finally
$$\eqalign{\scriptstyle \beta\,D_4 &\scriptstyle= b_4\,D_4 - D_1\,b_1 + D_2\,b_2 - D_3\,b_3 \cr 
&\scriptstyle=
   -h^{(0)}_{1}\,h^{(1)}_{1}\,h^{(2)}_{1}\,h^{(3)}_{1} -
    h^{(0)}_{1}\,h^{(1)}_{1}\,h^{(3)}_{1}\,h^{(2)}_{2} -
    h^{(0)}_{1}\,h^{(1)}_{1}\,h^{(2)}_{2}\,h^{(3)}_{2} -
    h^{(0)}_{1}\,h^{(1)}_{1}\,h^{(2)}_{2}\,h^{(3)}_{3} -
    h^{(0)}_{1}\,h^{(1)}_{1}\,h^{(2)}_{2}\,h^{(3)}_{4} \cr&\,\,\scriptstyle-
    h^{(0)}_{1}\,h^{(2)}_{1}\,h^{(3)}_{1}\,h^{(1)}_{2} -
    h^{(0)}_{1}\,h^{(2)}_{1}\,h^{(3)}_{1}\,h^{(1)}_{3} -
    h^{(0)}_{1}\,h^{(2)}_{1}\,h^{(3)}_{2}\,h^{(1)}_{3} -
    h^{(0)}_{1}\,h^{(2)}_{1}\,h^{(1)}_{3}\,h^{(3)}_{3} -
    h^{(0)}_{1}\,h^{(3)}_{1}\,h^{(1)}_{2}\,h^{(2)}_{2} \cr&\,\,\scriptstyle-
    h^{(0)}_{1}\,h^{(3)}_{1}\,h^{(1)}_{2}\,h^{(2)}_{3} -
    h^{(0)}_{1}\,h^{(3)}_{1}\,h^{(1)}_{2}\,h^{(2)}_{4} -
    h^{(0)}_{1}\,h^{(3)}_{1}\,h^{(2)}_{2}\,h^{(1)}_{3} -
    h^{(0)}_{1}\,h^{(3)}_{1}\,h^{(1)}_{3}\,h^{(2)}_{3} -
    h^{(0)}_{1}\,h^{(3)}_{1}\,h^{(1)}_{3}\,h^{(2)}_{4} \cr&\,\,\scriptstyle-
    h^{(0)}_{1}\,h^{(1)}_{2}\,h^{(2)}_{2}\,h^{(3)}_{4} -
    h^{(0)}_{1}\,h^{(2)}_{2}\,h^{(3)}_{2}\,h^{(1)}_{3} -
    h^{(0)}_{1}\,h^{(2)}_{2}\,h^{(1)}_{3}\,h^{(3)}_{3} -
    h^{(0)}_{1}\,h^{(2)}_{2}\,h^{(1)}_{3}\,h^{(3)}_{4} -
    h^{(0)}_{1}\,h^{(3)}_{2}\,h^{(1)}_{3}\,h^{(2)}_{3} \cr&\,\,\scriptstyle-
    h^{(0)}_{1}\,h^{(3)}_{2}\,h^{(1)}_{3}\,h^{(2)}_{4} -
    h^{(0)}_{1}\,h^{(1)}_{3}\,h^{(2)}_{3}\,h^{(3)}_{3} -
    h^{(0)}_{1}\,h^{(1)}_{3}\,h^{(2)}_{3}\,h^{(3)}_{4} -
    h^{(0)}_{1}\,h^{(1)}_{3}\,h^{(3)}_{3}\,h^{(2)}_{4} -
    h^{(0)}_{1}\,h^{(1)}_{3}\,h^{(2)}_{4}\,h^{(3)}_{4} \cr&\,\,\scriptstyle-
    h^{(1)}_{1}\,h^{(2)}_{1}\,h^{(3)}_{1}\,h^{(0)}_{2} -
    h^{(1)}_{1}\,h^{(2)}_{1}\,h^{(3)}_{1}\,h^{(0)}_{3} -
    h^{(1)}_{1}\,h^{(2)}_{1}\,h^{(3)}_{1}\,h^{(0)}_{4} -
    h^{(1)}_{1}\,h^{(2)}_{1}\,h^{(3)}_{2}\,h^{(0)}_{4} -
    h^{(1)}_{1}\,h^{(3)}_{1}\,h^{(0)}_{2}\,h^{(2)}_{2} \cr&\,\,\scriptstyle-
    h^{(1)}_{1}\,h^{(3)}_{1}\,h^{(2)}_{2}\,h^{(0)}_{3} -
    h^{(1)}_{1}\,h^{(3)}_{1}\,h^{(2)}_{2}\,h^{(0)}_{4} -
    h^{(1)}_{1}\,h^{(3)}_{1}\,h^{(0)}_{3}\,h^{(2)}_{3} -
    h^{(1)}_{1}\,h^{(3)}_{1}\,h^{(2)}_{3}\,h^{(0)}_{4} -
    h^{(1)}_{1}\,h^{(0)}_{2}\,h^{(2)}_{2}\,h^{(3)}_{2} \cr&\,\,\scriptstyle-
    h^{(1)}_{1}\,h^{(0)}_{2}\,h^{(2)}_{2}\,h^{(3)}_{3} -
    h^{(1)}_{1}\,h^{(0)}_{2}\,h^{(2)}_{2}\,h^{(3)}_{4} -
    h^{(1)}_{1}\,h^{(2)}_{2}\,h^{(3)}_{2}\,h^{(0)}_{4} -
    h^{(1)}_{1}\,h^{(2)}_{2}\,h^{(0)}_{3}\,h^{(3)}_{3} -
    h^{(1)}_{1}\,h^{(2)}_{2}\,h^{(0)}_{3}\,h^{(3)}_{4} \cr&\,\,\scriptstyle-
    h^{(1)}_{1}\,h^{(2)}_{2}\,h^{(3)}_{3}\,h^{(0)}_{4} -
    h^{(1)}_{1}\,h^{(2)}_{2}\,h^{(0)}_{4}\,h^{(3)}_{4} -
    h^{(1)}_{1}\,h^{(3)}_{2}\,h^{(2)}_{3}\,h^{(0)}_{4} -
    h^{(1)}_{1}\,h^{(2)}_{3}\,h^{(3)}_{3}\,h^{(0)}_{4} -
    h^{(1)}_{1}\,h^{(2)}_{3}\,h^{(0)}_{4}\,h^{(3)}_{4} \cr&\,\,\scriptstyle-
    h^{(2)}_{1}\,h^{(3)}_{1}\,h^{(0)}_{2}\,h^{(1)}_{2} -
    h^{(2)}_{1}\,h^{(3)}_{1}\,h^{(0)}_{2}\,h^{(1)}_{3} -
    h^{(2)}_{1}\,h^{(3)}_{1}\,h^{(0)}_{2}\,h^{(1)}_{4} -
    h^{(2)}_{1}\,h^{(3)}_{1}\,h^{(1)}_{2}\,h^{(0)}_{3} -
    h^{(2)}_{1}\,h^{(3)}_{1}\,h^{(1)}_{2}\,h^{(0)}_{4} \cr&\,\,\scriptstyle-
    h^{(2)}_{1}\,h^{(3)}_{1}\,h^{(0)}_{3}\,h^{(1)}_{3} -
    h^{(2)}_{1}\,h^{(3)}_{1}\,h^{(0)}_{3}\,h^{(1)}_{4} -
    h^{(2)}_{1}\,h^{(3)}_{1}\,h^{(1)}_{3}\,h^{(0)}_{4} -
    h^{(2)}_{1}\,h^{(3)}_{1}\,h^{(0)}_{4}\,h^{(1)}_{4} -
    h^{(2)}_{1}\,h^{(0)}_{2}\,h^{(1)}_{3}\,h^{(3)}_{3} \cr&\,\,\scriptstyle-
    h^{(2)}_{1}\,h^{(1)}_{2}\,h^{(3)}_{2}\,h^{(0)}_{4} -
    h^{(2)}_{1}\,h^{(3)}_{2}\,h^{(1)}_{3}\,h^{(0)}_{4} -
    h^{(2)}_{1}\,h^{(3)}_{2}\,h^{(0)}_{4}\,h^{(1)}_{4} -
    h^{(2)}_{1}\,h^{(0)}_{3}\,h^{(1)}_{3}\,h^{(3)}_{3} -
    h^{(2)}_{1}\,h^{(1)}_{3}\,h^{(3)}_{3}\,h^{(0)}_{4} \cr&\,\,\scriptstyle-
    h^{(2)}_{1}\,h^{(3)}_{3}\,h^{(0)}_{4}\,h^{(1)}_{4} -
    h^{(3)}_{1}\,h^{(0)}_{2}\,h^{(1)}_{2}\,h^{(2)}_{2} -
    h^{(3)}_{1}\,h^{(0)}_{2}\,h^{(1)}_{2}\,h^{(2)}_{4} -
    h^{(3)}_{1}\,h^{(0)}_{2}\,h^{(2)}_{2}\,h^{(1)}_{3} -
    h^{(3)}_{1}\,h^{(0)}_{2}\,h^{(2)}_{2}\,h^{(1)}_{4} \cr&\,\,\scriptstyle-
    h^{(3)}_{1}\,h^{(0)}_{2}\,h^{(1)}_{3}\,h^{(2)}_{4} -
    h^{(3)}_{1}\,h^{(1)}_{2}\,h^{(2)}_{2}\,h^{(0)}_{3} -
    h^{(3)}_{1}\,h^{(1)}_{2}\,h^{(2)}_{2}\,h^{(0)}_{4} -
    h^{(3)}_{1}\,h^{(1)}_{2}\,h^{(0)}_{3}\,h^{(2)}_{3} -
    h^{(3)}_{1}\,h^{(1)}_{2}\,h^{(0)}_{3}\,h^{(2)}_{4} \cr&\,\,\scriptstyle-
    h^{(3)}_{1}\,h^{(1)}_{2}\,h^{(2)}_{3}\,h^{(0)}_{4} -
    h^{(3)}_{1}\,h^{(1)}_{2}\,h^{(0)}_{4}\,h^{(2)}_{4} -
    h^{(3)}_{1}\,h^{(2)}_{2}\,h^{(0)}_{3}\,h^{(1)}_{3} -
    h^{(3)}_{1}\,h^{(2)}_{2}\,h^{(0)}_{3}\,h^{(1)}_{4} -
    h^{(3)}_{1}\,h^{(2)}_{2}\,h^{(1)}_{3}\,h^{(0)}_{4} \cr&\,\,\scriptstyle-
    h^{(3)}_{1}\,h^{(2)}_{2}\,h^{(0)}_{4}\,h^{(1)}_{4} -
    h^{(3)}_{1}\,h^{(0)}_{3}\,h^{(1)}_{3}\,h^{(2)}_{3} -
    h^{(3)}_{1}\,h^{(0)}_{3}\,h^{(1)}_{3}\,h^{(2)}_{4} -
    h^{(3)}_{1}\,h^{(0)}_{3}\,h^{(2)}_{3}\,h^{(1)}_{4} -
    h^{(3)}_{1}\,h^{(0)}_{3}\,h^{(1)}_{4}\,h^{(2)}_{4} \cr&\,\,\scriptstyle-
    h^{(3)}_{1}\,h^{(1)}_{3}\,h^{(2)}_{3}\,h^{(0)}_{4} -
    h^{(3)}_{1}\,h^{(1)}_{3}\,h^{(0)}_{4}\,h^{(2)}_{4} -
    h^{(3)}_{1}\,h^{(2)}_{3}\,h^{(0)}_{4}\,h^{(1)}_{4} -
    h^{(3)}_{1}\,h^{(0)}_{4}\,h^{(1)}_{4}\,h^{(2)}_{4} -
    h^{(0)}_{2}\,h^{(1)}_{2}\,h^{(2)}_{2}\,h^{(3)}_{2} \cr&\,\,\scriptstyle-
    h^{(0)}_{2}\,h^{(1)}_{2}\,h^{(2)}_{2}\,h^{(3)}_{4} -
    h^{(0)}_{2}\,h^{(2)}_{2}\,h^{(3)}_{2}\,h^{(1)}_{3} -
    h^{(0)}_{2}\,h^{(2)}_{2}\,h^{(3)}_{2}\,h^{(1)}_{4} -
    h^{(0)}_{2}\,h^{(2)}_{2}\,h^{(1)}_{3}\,h^{(3)}_{3} -
    h^{(0)}_{2}\,h^{(2)}_{2}\,h^{(1)}_{3}\,h^{(3)}_{4} \cr&\,\,\scriptstyle-
    h^{(0)}_{2}\,h^{(2)}_{2}\,h^{(3)}_{3}\,h^{(1)}_{4} -
    h^{(0)}_{2}\,h^{(2)}_{2}\,h^{(1)}_{4}\,h^{(3)}_{4} -
    h^{(0)}_{2}\,h^{(3)}_{2}\,h^{(1)}_{3}\,h^{(2)}_{4} -
    h^{(0)}_{2}\,h^{(1)}_{3}\,h^{(3)}_{3}\,h^{(2)}_{4} -
    h^{(0)}_{2}\,h^{(1)}_{3}\,h^{(2)}_{4}\,h^{(3)}_{4} \cr&\,\,\scriptstyle-
    h^{(1)}_{2}\,h^{(2)}_{2}\,h^{(3)}_{2}\,h^{(0)}_{4} -
    h^{(1)}_{2}\,h^{(2)}_{2}\,h^{(0)}_{3}\,h^{(3)}_{4} -
    h^{(1)}_{2}\,h^{(2)}_{2}\,h^{(0)}_{4}\,h^{(3)}_{4} -
    h^{(1)}_{2}\,h^{(3)}_{2}\,h^{(2)}_{3}\,h^{(0)}_{4} -
    h^{(1)}_{2}\,h^{(3)}_{2}\,h^{(0)}_{4}\,h^{(2)}_{4} \cr&\,\,\scriptstyle-
    h^{(1)}_{2}\,h^{(2)}_{3}\,h^{(0)}_{4}\,h^{(3)}_{4} -
    h^{(2)}_{2}\,h^{(3)}_{2}\,h^{(1)}_{3}\,h^{(0)}_{4} -
    h^{(2)}_{2}\,h^{(3)}_{2}\,h^{(0)}_{4}\,h^{(1)}_{4} -
    h^{(2)}_{2}\,h^{(0)}_{3}\,h^{(1)}_{3}\,h^{(3)}_{3} -
    h^{(2)}_{2}\,h^{(0)}_{3}\,h^{(1)}_{3}\,h^{(3)}_{4} \cr&\,\,\scriptstyle-
    h^{(2)}_{2}\,h^{(0)}_{3}\,h^{(1)}_{4}\,h^{(3)}_{4} -
    h^{(2)}_{2}\,h^{(1)}_{3}\,h^{(3)}_{3}\,h^{(0)}_{4} -
    h^{(2)}_{2}\,h^{(1)}_{3}\,h^{(0)}_{4}\,h^{(3)}_{4} -
    h^{(2)}_{2}\,h^{(3)}_{3}\,h^{(0)}_{4}\,h^{(1)}_{4} -
    h^{(2)}_{2}\,h^{(0)}_{4}\,h^{(1)}_{4}\,h^{(3)}_{4} \cr&\,\,\scriptstyle-
    h^{(3)}_{2}\,h^{(1)}_{3}\,h^{(2)}_{3}\,h^{(0)}_{4} -
    h^{(3)}_{2}\,h^{(1)}_{3}\,h^{(0)}_{4}\,h^{(2)}_{4} -
    h^{(3)}_{2}\,h^{(2)}_{3}\,h^{(0)}_{4}\,h^{(1)}_{4} -
    h^{(3)}_{2}\,h^{(0)}_{4}\,h^{(1)}_{4}\,h^{(2)}_{4} -
    h^{(0)}_{3}\,h^{(1)}_{3}\,h^{(2)}_{3}\,h^{(3)}_{3} \cr&\,\,\scriptstyle-
    h^{(0)}_{3}\,h^{(1)}_{3}\,h^{(3)}_{3}\,h^{(2)}_{4} -
    h^{(0)}_{3}\,h^{(1)}_{3}\,h^{(2)}_{4}\,h^{(3)}_{4} -
    h^{(1)}_{3}\,h^{(2)}_{3}\,h^{(3)}_{3}\,h^{(0)}_{4} -
    h^{(1)}_{3}\,h^{(2)}_{3}\,h^{(0)}_{4}\,h^{(3)}_{4} -
    h^{(1)}_{3}\,h^{(3)}_{3}\,h^{(0)}_{4}\,h^{(2)}_{4} \cr&\,\,\scriptstyle-
    h^{(1)}_{3}\,h^{(0)}_{4}\,h^{(2)}_{4}\,h^{(3)}_{4} -
    h^{(2)}_{3}\,h^{(3)}_{3}\,h^{(0)}_{4}\,h^{(1)}_{4} -
    h^{(2)}_{3}\,h^{(0)}_{4}\,h^{(1)}_{4}\,h^{(3)}_{4} -
    h^{(3)}_{3}\,h^{(0)}_{4}\,h^{(1)}_{4}\,h^{(2)}_{4} -
    h^{(0)}_{4}\,h^{(1)}_{4}\,h^{(2)}_{4}\,h^{(3)}_{4}
\cr}$$
and note that all 125 terms occur with a minus sign.  Since the $h^{(j)}_{i}$ are
all positive, this implies that $\beta\,D_4<0$, and since $D_4>0$, it follows that
$\beta<0$.

\definition{Final Note} If a single monomial term in any of
the expansions for $D_2$ or $D_4$ or $\beta\,D_4$ in terms of the
$h_i^{(j)}$'s had a negative sign, or if a single monomial term in any
of the expansions for $D_1$ or $D_3$ had a single negative term it
would be possible to choose a set of $h_i^{(j)}$'s so that at least
one of the $(-1)^k\,D_k$'s or $\beta$ would be positive (and hence that
$g_4$ would also be positive), which would lead to a realization
of the permutation set $Q_0$.
\enddefinition

\enddocument
\end